\documentclass[reqno]{amsart}

\usepackage{amsmath,amssymb}

\numberwithin{equation}{section}
\theoremstyle{plain}
\newtheorem*{thm}{Theorem}
\newtheorem*{cor}{Corollary}
\newtheorem*{lem}{Lemma}
\newtheorem*{prop}{Proposition}
\newtheorem*{Def}{Definition}
\theoremstyle{remark}

\begin{document}
\title{Askey-Wilson polynomials: an affine Hecke algebraic approach}

\date{January 6th, 2000}

\author{Masatoshi Noumi}
\address{Masatoshi Noumi, Department of Mathematics, Faculty of Science,
Kobe
University, Rokko, Kobe 657, Japan.}
\email{noumi@math.s.kobe-u.ac.jp}

\author{Jasper V. Stokman}
\address{Jasper V. Stokman, Centre de math{\'e}matiques de Jussieu,
Universit{\'e}
Paris 6 Pierre et Marie Curie, 4, Place Jussieu, F-75252 Paris Cedex
05, France.}
\email{stokman@math.jussieu.fr}

%\subjclass 2000: Primary 33D45, 33D80; Secondary 33D52


%%%%%%%%%%%%%%%%%%%%%%%%%%%%%%%%%%%%%%%%%%%%%%%%%%%%%%%%%%%%%%%%%
%%                                                             %%
%%                        Abstract                             %%
%%                                                             %%
%%%%%%%%%%%%%%%%%%%%%%%%%%%%%%%%%%%%%%%%%%%%%%%%%%%%%%%%%%%%%%%%%

\begin{abstract}
We study Askey-Wilson type polynomials using representation theory of
the double affine Hecke al\-ge\-bra. In par\-ti\-cu\-lar, we prove
bi-or\-tho\-go\-na\-li\-ty relations for non-sym\-me\-tric and
anti-symmetric Askey-Wilson polynomials with respect to a
complex measure. We give duality properties of the
non-symmetric Askey-Wilson polynomials, and we show how the
non-symmetric Askey-Wilson polynomials can be
created from Sahi's intertwiners.
The diagonal terms associated to the bi-orthogonality relations
(which replace the notion of quadratic norm evaluations for
orthogonal polynomials) are expressed in terms of residues of
the complex weight function using intertwining properties of
the non-symmetric Askey-Wilson transform under the action of the
double affine Hecke algebra. We evaluate the constant term, which is
essentially the well-known Askey-Wilson integral, using shift
operators. We furthermore show how these results
reduce to well-known properties of the symmetric Askey-Wilson
polynomials, as were originally derived by Askey and Wilson using
basic hypergeometric series theory.
\end{abstract}

%%%%%%%%%%%%%%%%%%%%%%%%%%%%%%%%%%%%%%%%%%%%%%%%%%%%%%%%%%%%%%%%%%%

\maketitle

%%%%%%%%%%%%%%%%%%%%%%%%%%%%%%%%%%%%%%%%%%%%%%%%%%%%%%%%%%%%%%%%
%%                                                            %%
%%                    Introduction                            %%
%%                                                            %%
%%%%%%%%%%%%%%%%%%%%%%%%%%%%%%%%%%%%%%%%%%%%%%%%%%%%%%%%%%%%%%%%

\section{Introduction}\label{section1}

\subsection {}
%1.1

Due to work of Cherednik \cite{C1}--\cite{C5}, Macdonald \cite{M2},
Noumi \cite{N} and Sahi \cite{S}, one can associate to
every irreducible affine root system certain families of orthogonal
polynomials (all closely related to the Macdonald polynomials), and
prove their basic properties using a fundamental representation of
the affine Hecke algebra in terms of difference-reflection
operators. In this paper, we consider a rank one example of this
theory in detail. The example is connected with a rank one
non-reduced irreducible affine root system which has four orbits
under the action of the associated affine Weyl group. The family of
symmetric orthogonal polynomials
associated to this particular affine root system
is the celebrated four parameter family
of Askey-Wilson polynomials, see \cite{AW}.

\subsection {}
%1.2

The four parameter family of Askey-Wilson polynomials has played
an important and central role in the theory of basic hypergeometric
orthogonal polynomials.
In fact, up to date they seem to be the most general family
of basic hypergeometric orthogonal polynomials which satisfy the
additional requirement that they are joint eigenfunctions of
a second-order $q$-difference operator.
We use the link between the Askey-Wilson polynomials and the most general
non-reduced affine root system of rank one (see 1.1)
to derive in this paper the basic properties
of the Askey-Wilson polynomials (and more!) from the algebraic
structure of the associated (double) affine Hecke algebra.

\subsection {}

%1.3
We introduce a Cherednik-Dunkl type difference-reflection operator $Y$
using the fundamental representation of the affine Hecke algebra of
type $\widetilde{A}_1$. The fundamental representation
was defined by Noumi \cite{N} in the higher rank case, see also Sahi
\cite{S}. Sahi's \cite{S} non-symmetric Askey-Wilson polynomials are
defined as the eigenfunctions of the Cherednik-Dunkl
operator $Y$. They form a linear basis of the Laurent
polynomials in one variable. We explicitly indicate their connection
with the symmetric Askey-Wilson polynomials as
originally defined by Askey and Wilson in \cite{AW}. In particular,
we give explicit expressions for the non-symmetric
Askey-Wilson polynomials as a sum of two terminating
balanced ${}_4\phi_3$'s (here
${}_r\phi_s$ is the basic hypergeometric series,
see Gasper and Rahman \cite{GR}
for the definition). All the other results in this paper
are derived without using the explicit series expansions of the
(non-)symmetric
Askey-Wilson polynomials in terms of basic hypergeometric series.

\subsection {}

%1.4
We derive bi-orthogonality relations for the non-symmetric
Askey-Wilson polynomials by explicitly computing the adjoint of the
Cherednik-Dunkl operator $Y$ with respect to an explicit, complex measure.
By a kind of symmetrization procedure, the bi-orthogonality
relations imply Askey and Wilson's \cite{AW}
orthogonality relations for the symmetric
Askey-Wilson polynomials.

\subsection {}

%1.5
We shortly describe Sahi's
\cite{S} results how an anti-isomorphism of
the double affine Hecke algebra gives rise to a duality between
the spectral and the
geometric parameter of the (non-)symmetric Askey-Wilson polynomial.
We use the duality to determine explicit intertwining properties
of the action of the double affine Hecke algebra under the
non-symmetric Askey-Wilson transform. Here the non-symmetric
Askey-Wilson transform is a generalized Fourier transform
which is defined with respect to the bi-orthogonality measure of the
non-symmetric Askey-Wilson polynomials. This leads to the
evaluation of the diagonal terms of
the (bi-)orthogonality relations in terms of certain
residues of the complex weight function.
This technique is motivated by Cherednik's \cite{C4}
approach to the study of the diagonal terms for non-symmetric
Macdonald polynomials, in which he rewrites part of the action of
the double affine Hecke algebra on non-symmetric Macdonald polynomials
in terms of explicit operators acting on the spectral parameter.

\subsection {}

%1.6
To complete the explicit computation of the diagonal terms,
we need the evaluation of the constant term and the evaluation
of the non-symmetric Askey-Wilson polynomial in a specific point
(the latter playing a fundamental role in the duality arguments).
We evaluate the constant term, which is essentially the
well-known Askey-Wilson integral (see \cite{AW}), using shift operators.
For the evaluation of the non-symmetric Askey-Wilson polynomial in a
specific point we use a Rodrigues type formula for the non-symmetric
Askey-Wilson polynomial in terms of Sahi's \cite{S} intertwiners.

\subsection {}

%1.7
The purpose of this paper is two-fold. First of all, we would like
to show the power of the Cherednik-Macdonald theory in the study of
basic hypergeometric orthogonal polynomials.
It not only shows that all the basic properties of the Askey-Wilson
polynomials can be obtained
by natural algebraic manipulations, but it also reveals new and
important insights in the structure of the
Askey-Wilson polynomials.

Secondly, the affine Hecke algebraic
approach also works for multivariable Askey-Wilson polynomials,
the so called Koornwinder polynomials \cite{Ko1} (which are associated with
a higher rank non-reduced affine root system).
The structure of the proofs in the higher rank setting are essentially
the same, although the technicalities are more involved.
Only part of the Cherednik-Macdonald theory
associated with non-reduced affine root systems has been written
down explicitly at this moment, see Noumi \cite{N} and Sahi \cite{S}.
In our opinion, this paper can serve as one of the building blocks
for obtaining a full understanding of the Cherednik-Macdonald theory
in the case of non-reduced affine root systems.
The higher rank case will be treated in an upcoming paper of
the second author.

\subsection {}

%1.8
In view of our aims described in 1.7, we have chosen to make this
paper fairly self-contained.
In particular, we have included some of the proofs
of Noumi \cite{N} and Sahi \cite{S}, restricted to our present rank one
setting.
Furthermore, we have included several proofs which are
fairly straightforward modifications from Cherednik's \cite{C1}--\cite{C5}
and Macdonald's \cite{M2} work in case of Macdonald polynomials, see
also for instance Opdam's \cite{O} lecture notes for the classical $q=1$
setting.

\subsection {}

%1.9
Finally we would like to point out the close connection with
the paper  \cite{KM} of Kalnins and Miller. In \cite{KM} the
Askey-Wilson second order $q$-difference operator is written as the
composition of a first order $q$-difference operator with its adjoint.
This decomposition leads naturally to proofs of the orthogonality
relations and of the quadratic norm evaluations for the symmetric
Askey-Wilson polynomials (using shift principles).
In our paper we use similar techniques, but we decompose the
Askey-Wilson second order $q$-difference operator now as a sum
of a difference-reflection operator (the Cherednik-Dunkl operator $Y$)
and its inverse. This decomposition has the advantage that the
Cherednik-Dunkl operator $Y$ itself satisfies a
self-adjointness
property with respect to a suitable extension of the
Askey-Wilson orthogonality measure to a
complex measure space for {\it non-symmetric} functions.
This extra property of $Y$ naturally leads to the introduction of
non-symmetric and anti-symmetric Askey-Wilson
polynomials and to their bi-orthogonality relations.

\subsection {}

%1.9
{\it Notations:} We use Gasper and Rahman's \cite{GR} notations for basic
hypergeometric series and $q$-shifted factorials.
We write ${\mathbb{Z}}_+=\{0,1,2,\ldots\}$ for the positive
integers and ${\mathbb{N}}=\{1,2,\ldots\}$ for the strictly positive
integers.

\subsection {}

%1.10
{\it Acknowledgments:} The second author is supported by a
NWO-TALENT stipendium of the Netherlands Organization for
Scientific Research (NWO). Part of the research was done while the
second author was supported by the EC TMR network ``Algebraic Lie
Representations'', grant no. ERB FMRX-CT97-0100.

%%%%%%%%%%%%%%%%%%%%%%%%%%%%%%%%%%%%%%%%%%%%%%%%%%%%%%%%%%%%%%%%
%%                                                            %%
%%    The Dunkl-Cherednik difference-reflection operators     %%
%%                                                            %%
%%%%%%%%%%%%%%%%%%%%%%%%%%%%%%%%%%%%%%%%%%%%%%%%%%%%%%%%%%%%%%%%

\section{The Dunkl-Cherednik difference-reflection operators}

\subsection {}

%2.1
Let $\widehat{\mathbb{R}}$ be the
vector space consisting of affine, linear transformation from
${\mathbb{R}}$ to ${\mathbb{R}}$.
We identify $\widehat{\mathbb{R}}$ with
${\mathbb{R}}\oplus {\mathbb{R}}\delta$,
where $\delta$ is the function identically one, and where
${\mathbb{R}}$ acts by multiplication on itself.
We introduce and study in this section a particular example of a rank one
affine root system $S\subset \widehat{\mathbb{R}}$. See Macdonald
\cite{M1} for the general discussion of affine root systems.

\subsection {}

%2.2
Let $\langle .,. \rangle$
be the positive semi-definite form on $\widehat{\mathbb{R}}$ defined by
\[ \langle \lambda+\mu\delta, \lambda'+\mu'\delta\rangle=
\lambda\lambda',\qquad \lambda,\lambda',\mu,\mu'\in {\mathbb{R}}.
\]
Then we can associate to every $0\not=f=\lambda+\mu\delta\in
\widehat{\mathbb{R}}\setminus {\mathbb{R}}\delta$ the involution
$s_f: \widehat{\mathbb{R}}\rightarrow \widehat{\mathbb{R}}$ defined by
\[ s_f(g)=g-\langle g,f^\vee\rangle f, \qquad g\in
\widehat{\mathbb{R}},
\]
where $f^\vee:=2f/\langle f,f\rangle$.
Observe that $s_f$ is an isometry with respect to $\langle .,.
\rangle$. In fact, we have $s_f(g)=g\circ \widetilde{s}_f$, where
$\widetilde{s}_f\in \widehat{\mathbb{R}}$
is the reflection in $f^{-1}(0)$.

\subsection {}

%2.3
We define now a subset $S\subset \widehat{\mathbb{R}}$
by
\[ S=\{ \pm 1+\frac{m}{2}\delta,\,\,\, \pm 2 +m\delta \,\, |
\,\,
m\in {\mathbb{Z}} \},
\]
and we write ${\mathcal{W}}={\mathcal{W}}(S)$ for the subgroup
of invertible linear
transformations of $\widehat{\mathbb{R}}$ generated by $s_f$, $f\in S$.
Then $S\subset \widehat{\mathbb{R}}$ is an affine root system.
In particular, $\langle f,g^\vee\rangle\in {\mathbb{Z}}$
for all $f,g\in S$, and $S$ is stable under
the action of the affine Weyl group ${\mathcal{W}}(S)$.

\subsection {}

%2.4
The {\it gradient root system} $\Sigma$ of $S$ is the projection
of $S$ on ${\mathbb{R}}$ along the direct
sum decomposition
$\widehat{\mathbb{R}}={\mathbb{R}}\oplus {\mathbb{R}}\delta$.
Here $\Sigma=\{ \pm 1, \pm 2\}$,
which is a non-reduced root system of type $BC_1$, with associated
Weyl group $W=\{1,s_1\}=\{ \pm 1\}$. Observe that
$W\subset {\mathcal{W}}$ acts on $\widehat{\mathbb{R}}$ by
$(\pm 1)(\lambda+\mu\delta)=\pm\lambda+\mu\delta$.

\subsection {}

%2.5
Let $a_0=\delta-2\in S$ and
$a_1=2\in S$. Observe that $a_0^\vee=
a_0/2=\frac{1}{2}\delta-1\in S$, and
$a_1^\vee=a_1/2=1\in S$.
Then $\{a_0^\vee, a_1^\vee\}$ forms a basis of the affine root
system $S$.
We write $S^+$ for the positive roots in
$S$ with respect to $\{a_0^\vee, a_1^\vee\}$, so that $S=S^+\cup (-S^+)$
disjoint union. We furthermore set $\Sigma^+=\{a_1^\vee, a_1\}$, which are
the positive roots of $\Sigma$ with respect to the basis $\{a_1^\vee\}$.
Observe that $S^+=\Sigma^+\cup \{ f\in S \, | \, f(0)>0\}$.

\subsection {}

%2.6
The affine Weyl group ${\mathcal{W}}$ is generated
by the simple
reflections $s_0:=s_{a_0}=s_{a_0^\vee}$  and $s_1=s_{a_1}=s_{a_1^\vee}$,
while $W$ is generated by $s_1$.
Observe that
$s_1s_0=\tau(1)$, where $\tau(\mu)$ is the translation operator
$\tau(\mu)f=f+\langle \mu, f\rangle\delta$ for $f\in \widehat{\mathbb{R}}$
and $\mu\in {\mathbb{R}}$.
In particular, $s_1s_0$ has infinite order in ${\mathcal{W}}$ and
\[ {\mathcal{W}}=W\ltimes \tau({\mathbb{Z}}).
\]
Furthermore, ${\mathcal{W}}$ is isomorphic to the Coxeter group with
two
generators $s_0$, $s_1$ and relations $s_0^2=1$,
$s_1^2=1$.

\subsection {}

%2.7
The affine root system $S$ has four ${\mathcal{W}}$-orbits, namely
\begin{equation*}
\begin{split}
S_s^1&={\mathcal{W}}a_0^\vee=
\bigl(\frac{1}{2}+{\mathbb{Z}}\bigr)\delta\pm 1,\qquad
S_s^2={\mathcal{W}}a_1^\vee={\mathbb{Z}}\delta\pm 1,\\
S_l^1&={\mathcal{W}}a_0=\bigl(1+2{\mathbb{Z}}\bigr)\delta\pm 2,\qquad
S_l^2={\mathcal{W}}a_1=2{\mathbb{Z}}\delta\pm 2.
\end{split}
\end{equation*}

\subsection {}

%2.8
We have the disjoint union
$S=R^\vee\cup R$, with
$R=S_l^1\cup S_l^2$ a reduced, irreducible affine root system
with basis $\{ a_0, a_1\}$, affine Weyl group ${\mathcal{W}}$, and gradient
root system $\{ \pm a_1\}$, and with
$R^\vee=S_s^1\cup S_s^2$ the corresponding affine co-root system. The
co-root system $R^\vee$ is
a reduced, affine root system with basis $\{ a_0^\vee, a_1^\vee\}$,
affine Weyl group
${\mathcal{W}}$, and gradient root system $\{ \pm a_1^\vee\}$.
Similarly as for $S$, see 2.5, the fixed choice of basis
give rise to a decomposition of $R$ and $R^\vee$ in
positive
and negative roots (the positive roots are denoted by $R^+$ and
$R^{\vee,+}$, respectively).
For example, we have $R^+=\{a_1\}\cup \{\pm a_1+{\mathbb{N}}\delta\}$.

\subsection {}

%2.9
Let $\omega\in \widehat{\mathbb{R}}$ be the involution
$\omega(x):=\frac{1}{2}-x$ ($x\in {\mathbb{R}}$), and consider
$\omega$ as an involution of
$\widehat{\mathbb{R}}$ by $\omega: f\mapsto f\circ
\omega$
for $f\in \widehat{\mathbb{R}}$. Then $\omega$
preserves $S$. Furthermore,
$\omega=s_1\tau(\frac{1}{2})$, $\omega(a_i)=a_{1-i}$
and $\omega s_i\omega=s_{1-i}$ for $i=0,1$.
The subgroup ${\mathcal{W}}^e$ of the invertible linear transformations
of $\widehat{\mathbb{R}}$ generated
by ${\mathcal{W}}$ and $\omega$ is called the
{\it extended affine Weyl group}.
It is isomorphic to ${\mathcal{W}}\rtimes \Omega$, where $\Omega$ is
the subgroup of order two generated by $\omega$.

\subsection {}

%2.10
Set ${\mathcal{A}}:= {\mathbb{C}}[x^{\pm 1}]$ for the algebra of
Laurent polynomials in one indeterminate $x$.
We set $x^f:=q^{\lambda}x^{\mu}\in {\mathcal{A}}$
for $f=\mu+\lambda\delta\in
{\mathbb{Z}}+{\mathbb{R}}\delta$,
where $q$ is a fixed non-zero complex number.
Observe that $x^a\in {\mathcal{A}}$ is well defined for $a\in S$ and that
${\mathcal{W}}^e$ preserves ${\mathbb{Z}}+{\mathbb{R}}\delta$.
Furthermore, $w(x^{\mu}):=x^{w(\mu)}$ for $\mu\in {\mathbb{Z}}$ and $w\in
{\mathcal{W}}^e$
extends to an action of ${\mathcal{W}}^e$ on ${\mathcal{A}}$ by linearity.
In particular,
\[ s_0(x^m)=q^{m}x^{-m},\qquad s_1(x^{m})=x^{-m},\qquad m\in {\mathbb{Z}}.
\]
Observe that $\tau(\mu)$ ($\mu\in {\mathbb{Z}}$) acts as a
$q$-difference operator: $\tau(\mu)(x^m)=q^{\mu m}x^m$ for all $m\in
{\mathbb{Z}}$.

\subsection {}

%2.11
A {\it multiplicity function} $\underline{t}=\{t_a\}_{a\in S}$ of
$S$
is a choice of non-zero complex numbers $t_a$ ($a\in S$)
such that $t_{w(a)}=t_a$
for all $a\in S$ and all $w\in W$.
We use the convention that $t_f=1$ for all $f\in
\widehat{\mathbb{R}}\setminus
S$.
A multiplicity function $\underline{t}$ of $S$ is determined by the values
$k_0:=t_{a_0}$, $u_0:=t_{a_0^\vee}$, $k_1:=t_{a_1}$ and
$u_1:=t_{a_1^\vee}$, see 2.7.
Later on it will be necessary to impose (generic) conditions on the
parameters $t_f$ ($f\in S$) and on the deformation parameter $q$.
Until section 6 it suffices to assume that $|q|\not=1$ and that
$k_0^2, k_1^2, u_1^2\not\in \pm q^{\mathbb{Z}}$. These conditions are
assumed to hold throughout the remainder of the paper, unless
specified explicitly otherwise.

\subsection {}

%2.12
The {\it Hecke algebra} $H_0=H_0(k_1)$ of type $A_1$
is the unital, associative ${\mathbb{C}}$-algebra
with generator $T_1$ and relation $(T_1-k_1)(T_1+k_1^{-1})=0$.
Observe that $\{1, T_1\}$ is a
linear basis of $H_0$ and that $T_1$ is invertible in $H_0$
with inverse $T_1^{-1}=T_1+k_1^{-1}-k_1$.

\subsection {}

%2.13
The {\it affine Hecke algebra} $H=H(R;k_0,k_1)$ of type $\widetilde{A}_1$
is
the unital ${\mathbb{C}}$-algebra with generators $T_0$ and $T_1$
and relations
\[(T_i-k_i)(T_i+k_i^{-1})=0,\qquad i=0,1.
\]
Similarly as for $H_0$ we have that $T_i$ is invertible in $H$ with
inverse $T_i^{-1}=T_i+k_i^{-1}-k_i$.

\subsection {}

%2.14
For $w\in {\mathcal{W}}$, let $w=s_{i_1}s_{i_2}\cdots s_{i_r}$ be a
reduced
expression, i.e. a minimal expression of $w$ as product
of the simple
reflections $s_0$ and $s_1$. Then $T_w:=T_{i_1}T_{i_2}\cdots
T_{i_r}$ is well-defined and
$\{ T_w\}_{w\in {\mathcal{W}}}$ is a linear basis of
$H(R;k_0,k_1)$, see \cite{L}. In particular, we may regard $H_0$ as
a subalgebra of $H$.

\subsection {}

%2.15
We set
$Y:=T_{\tau(1)}=T_1T_0\in H$. By \cite{L}, we known that
$Y$ is algebraically independent in $H$. Let
${\mathbb{C}}[Y^{\pm 1}]\subset H$
be the commutative subalgebra generated by $Y^{\pm 1}$. Then
\[ H_0(k_1)\otimes {\mathbb{C}}[Y^{\pm 1}]\simeq
H(R;k_0,k_1) \simeq {\mathbb{C}}[Y^{\pm 1}]\otimes H_0(k_1)
\]
as linear spaces, where the isomorphisms are given by multiplication.
In particular, $\{ Y^{m}, Y^nT_1\}_{m,n\in {\mathbb{Z}}}$ and
$\{ Y^m, T_1Y^n\}_{m,n\in {\mathbb{Z}}}$ are linear bases of
$H(R;k_0,k_1)$, see \cite[proposition 3.7]{L}.

In the remainder of the paper we identify ${\mathcal{A}}$ with
${\mathbb{C}}[Y^{\pm 1}]$ as algebra by identifying the
indeterminate $x$ of ${\mathcal{A}}$ with $Y$.
In particular, we write $f(Y)=\sum_kc_kY^k\in {\mathbb{C}}[Y^{\pm 1}]$
for $f(x)=\sum_kc_kx^k\in {\mathcal{A}}$.

\subsection {}

%2.16
By Lusztig \cite[proposition 3.6]{L},
we have the fundamental commutation relations
\[T_1f(Y)-f(Y^{-1})T_1=\bigl((k_1-k_1^{-1})Y^2+(k_0-k_0^{-1})Y\bigr)
\left(\frac{f(Y^{-1})-f(Y)}{1-Y^2}\right)
\]
in $H(R;k_0,k_1)$ for all $f(Y)\in {\mathbb{C}}[Y^{\pm 1}]$.
Indeed, observe that if the formula holds for $f(Y)$ and $g(Y)$,
then it also holds for $f(Y)g(Y)$.
It thus suffices to prove it for $f(Y)=Y^{\pm 1}$,
in which case it follows from an
elementary computation using the definition of $Y$
and the quadratic relations for the $T_i$, see 2.15 and 2.13,
respectively.

\subsection {}

%2.17
The following result was observed by Sahi \cite[theorem 5.1]{S} in the
higher rank setting.
\begin{cor}[The non-affine intertwiner]
Set $S_1=[T_1,Y]=T_1Y-YT_1\in H$. Then
$f(Y)S_1=S_1f(Y^{-1})$ in $H$ for all $f(Y)\in {\mathbb{C}}[Y^{\pm 1}]$.
\end{cor}
\begin{proof}
This follows immediately from the definition of $S_1$
and from Lusztig's commutation relation 2.16.
\end{proof}

\subsection {}

%2.18
Another important consequence of Lusztig's commutation relation 2.16
is the following result.
\begin{cor}The affine Hecke algebra $H=H(R;k_0,k_1)$ acts on
${\mathbb{C}}[Y^{\pm 1}]$ by
\begin{equation*}
\begin{split}
T_1.g(Y)&=k_1g(Y^{-1})+\bigl((k_1-k_1^{-1})Y^2+(k_0-k_0^{-1})Y\bigr)
\left(\frac{g(Y^{-1})-g(Y)}{1-Y^2}\right)\\&=k_1g(Y)+k_1^{-1}
\frac{(1-k_0k_1Y^{-1})(1+k_0^{-1}k_1Y^{-1})}{(1-Y^{-2})}(g(Y^{-1})-g(Y)),\\
f(Y).g(Y)&=f(Y)g(Y)
\end{split}
\end{equation*}
for all $f(Y), g(Y)\in {\mathbb{C}}[Y^{\pm 1}]$.
\end{cor}
\begin{proof}
Let $\chi$ be the character of $H_0(k_1)$ which maps $T_1$ to
$k_1$.
By 2.15 we may identify the representation space of
the induced representation ${\hbox{Ind}}_{H_0}^{H}(\chi)=
H\otimes_{\chi}{\mathbb{C}}$ with ${\mathbb{C}}[Y^{\pm 1}]$.
By 2.16 the corresponding induced action
of $H$ on ${\mathbb{C}}[Y^{\pm 1}]$ is as indicated in the statement
of the corollary.
\end{proof}

\subsection {}

%2.19
We define linear operators
$\widehat{T}_i\in \hbox{End}_{\mathbb{C}}({\mathcal{A}})$ by
\begin{equation*}
\begin{split}
\widehat{T}_i:=&k_i+k_i^{-1}\frac{\bigl(1-k_iu_ix^{a_i^\vee}\bigr)
\bigl(1+k_iu_i^{-1}x^{a_i^\vee}\bigr)}{1-x^{a_i}}(s_i-\hbox{id})\\
=&k_is_i+\frac{(k_i-k_i^{-1})+(u_i-u_i^{-1})x^{a_i^\vee}}{(1-x^{a_i})}
(\hbox{id}-s_i), \qquad i=0,1.
\end{split}
\end{equation*}
The following theorem was proved by Noumi \cite[section 3]{N} in
the higher rank setting, see also \cite[section 2.3]{S}.
\begin{thm}
The application $T_i\mapsto \widehat{T}_i$ \textup{(}$i=0,1$\textup{)}
extends uniquely
to an algebra homomorphism $\pi_{\underline{t},q}:
H(R;k_0,k_1)\rightarrow
\hbox{End}_{{\mathbb{C}}}({\mathcal{A}})$.
\end{thm}
\begin{proof}
We identify ${\mathbb{C}}[Y^{\pm 1}]\subset H(R;u_1,k_1)$
with ${\mathcal{A}}$ as algebra by identifying $Y$ with $x^{-1}$.
Then it follows from corollary 2.18 (applied to $H(R;u_1,k_1)$)
that $(\widehat{T}_1-k_1)(\widehat{T}_1+k_1^{-1})=0$
in $\hbox{End}_{\mathbb{C}}({\mathcal{A}})$.
Conjugating $\widehat{T}_1$ with the involution $\omega$, see 2.9,
and replacing $k_1$ and $u_1$ by $k_0$ and $u_0$ respectively, we see
that $(\widehat{T}_0-k_0)(\widehat{T}_0+k_0^{-1})=0$ in
$\hbox{End}({\mathcal{A}})$. The theorem
follows, since we have shown that all
the defining relations 2.13 of $H(R;k_0,k_1)$ are satisfied
by the linear operators $\widehat{T}_i\in
\hbox{End}_{\mathbb{C}}({\mathcal{A}})$
($i=0,1$).
\end{proof}

\subsection {}

%2.20
Observe that the linear operator $\widehat{T}_0$ on ${\mathcal{A}}$
has a reflection and a $q$-difference part, while $\widehat{T}_1$
has only a reflection part, see  2.10.
The operators $\widehat{T}_i\in \hbox{End}_{\mathbb{C}}({\mathcal{A}})$
($i=0,1$) are called the {\it difference-reflection
operators associated with $S$}. In the remainder of
the paper, we simply write $T_i$ for the difference-reflection
operators $\widehat{T}_i$ ($i=0,1$) if no confusion is possible.
The operator
$Y=T_1T_0\in \hbox{End}_{\mathbb{C}}({\mathcal{A}})$ is called the
{\it Cherednik-Dunkl operator} associated with $S$.

\subsection {}

%2.21
The representation $\pi_{\underline{t},q}$
has two extra degrees of freedom $u_0$ and $u_1$ besides the
deformation parameter $q$
(which already lives on the affine Weyl group level,
see 2.10). The motivation to label these two
degrees of freedom in this particular way comes
from the theory of double affine Hecke algebras.
The {\it double affine Hecke algebra}
$\mathcal{H}(S;\underline{t};q)$
associated with the affine root system $S$ (see \cite{S})
is the subalgebra of $\hbox{End}_{\mathbb{C}}({\mathcal{A}})$
generated by $\pi_{\underline{t},q}(H(R;k_0,k_1))$ and
${\mathcal{A}}$,
where we consider ${\mathcal{A}}$ as a subalgebra of
$\hbox{End}_{\mathbb{C}}({\mathcal{A}})$
via its regular representation. We write
$f(z)\in \hbox{End}\bigl({\mathcal{A}}\bigr)$ for the Laurent
polynomial $f(x)\in {\mathcal{A}}$ regarded as a linear
endomorphism of ${\mathcal{A}}$.
By the second formula for the difference-reflection operators $T_i$
in 2.19 we have that
\[
f(z)T_i-T_i(s_if)(z)=\frac{(k_i-k_i^{-1})+(u_i-u_i^{-1})z^{a_i^\vee}}
{(1-z^{a_i})}\bigl(f(z)-(s_if)(z)\bigr),\qquad f\in {\mathcal{A}}
\]
in ${\mathcal{H}}(S;\underline{t},q)$ for $i=0,1$.

\subsection {}

%2.22
The labeling of the extra degrees of
freedom in the representation $\pi_{\underline{t},q}$
is now justified by the following theorem, together with 2.11.
\begin{thm}
$\mathcal{H}(S;\underline{t};q)$ is isomorphic as algebra
to the unital, associative  ${\mathbb{C}}$-algebra
${\mathcal{F}}(\underline{t};q)$ with generators $V_0^\vee, V_0,
V_1, V_1^\vee$ and relations:
\begin{enumerate}
\item{} The application $T_i\mapsto V_i$ for $i=0,1$
extends to an algebra homomorphism $H(R;k_0,k_1)\rightarrow
{\mathcal{F}}(\underline{t};q)$.
\item{} The application $T_i\mapsto V_i^\vee$ for $i=0,1$ extends
to
an algebra homomorphism $H(R;u_0,u_1)\rightarrow
{\mathcal{F}}(\underline{t};q)$.
\item{} \textup{(}Compatibility\textup{)}. $V_1^\vee V_1
V_0V_0^\vee=
q^{-1/2}$.
\end{enumerate}
The isomorphism $\phi: {\mathcal{F}}(\underline{t};q)\rightarrow
{\mathcal{H}}(S;\underline{t};q)$ is explicitly given by $\phi(V_i)= T_i$
\textup{(}$i=0,1$\textup{)},
$\phi(V_0^\vee)=T_0^{-1}z^{-a_0^\vee}=q^{-1/2}T_0^{-1}z$
and $\phi(V_1^\vee)=z^{-a_1^\vee}T_1^{-1}=z^{-1}T_1^{-1}$.
\end{thm}
The existence of the algebra homomorphism $\phi$ follows
by direct computations using 2.21. It is immediate that $\phi$
is surjective. The injectivity of $\phi$ requires a detailed study of the
difference-reflection operators $T_i$ associated with $S$. We give the
proof in
8.3.

%%%%%%%%%%%%%%%%%%%%%%%%%%%%%%%%%%%%%%%%%%%%%%%%%%%%%%%%%%%%
%%                                                        %%
%%                                                        %%
%%           Non-symmetric Askey-Wilson polynomials       %%
%%                                                        %%
%%                                                        %%
%%%%%%%%%%%%%%%%%%%%%%%%%%%%%%%%%%%%%%%%%%%%%%%%%%%%%%%%%%%%

\section{Non-symmetric Askey-Wilson polynomials}

\subsection {}

%3.1
For $a\in R$, let ${\mathcal{R}}(a)\in
\hbox{End}({\mathcal{A}})$
be the difference-reflection operator defined by
\[
{\mathcal{R}}(a):=t_as_a+t_a^{-1}\frac{(1-t_at_{a/2}x^{a/2})
(1+t_at_{a/2}^{-1}x^{a/2})}{(1-x^a)}\bigl(1-s_a\bigr).
\]
Then it is immediate that ${\mathcal{R}}(a_i)=T_is_i$ for
$i=0,1$, where $T_i$ is the
difference-reflection operator associated with $S$ (see 2.19).
Furthermore, we have $w{\mathcal{R}}(a)w^{-1}={\mathcal{R}}(w(a))$ for all
$w\in {\mathcal{W}}$ and all $a\in R$.
Since any $a\in R$ is conjugate to $a_0$ or $a_1$ under the
action of ${\mathcal{W}}$, we obtain from the quadratic relations
for the difference-reflection operators $T_i$ (see 2.13 and 2.19) that
\[ {\mathcal{R}}(a)^{-1}={\mathcal{R}}(-a)+(t_a-t_a^{-1})s_a,\qquad
a\in R.
\]

\subsection {}

%3.2
We define  a total order $\preceq$
on the basis of monomials $\{x^m\}_{m\in {\mathbb{Z}}}$ of ${\mathcal{A}}$
by
\[1\prec x^{-1}\prec x\prec x^{-2}\prec x^2\prec\cdots.
\]
Observe that this order is not well behaved under multiplication of
the monomials: if $x^{m_i}\preceq x^{n_i}$ ($i=1,2$), then not necessarily
$x^{m_1+m_2}\preceq x^{n_1+n_2}$.

\subsection {}

%3.3
Let $\epsilon: {\mathbb{Z}}\rightarrow \{\pm 1\}$ be the function
which sends a positive integer to $1$ and a strictly negative integer
to $-1$.
\begin{lem}
Let $a\in R$ be of the form $a=2+k\delta$ with $k\in {\mathbb{Z}}$
\textup{(}see 2.8\textup{)} and let $m\in {\mathbb{Z}}$. Then
\[{\mathcal{R}}(a)\bigl(x^m\bigr)=t_a^{\epsilon(m)}x^m+
\, \hbox{\textup{lower order terms w.r.t.}}\, \preceq.
\]
\end{lem}
\begin{proof}
For $a\in R$, let $D_a\in \hbox{End}({\mathcal{A}})$
be the divided difference operator defined by
\[ D_af:=\frac{f-s_af}{1-x^a}, \qquad f\in {\mathcal{A}}.
\]
Then $D_a(1)=0$ and
\begin{equation*}
\begin{split}
D_a\bigl(x^m\bigr)=
\begin{cases}
-x^{m-a}-x^{m-2a}-\cdots -x^{m-\langle m,a^\vee\rangle a}\qquad
&\hbox{ if } \langle m, a^\vee\rangle\in {\mathbb{N}},\\
x^m+x^{m+a}+\cdots +x^{m-(1+\langle m,a^\vee\rangle)a}\qquad
&\hbox{ if } \langle m, a^\vee\rangle\in -{\mathbb{N}}.
\end{cases}
\end{split}
\end{equation*}
Observe that $\langle m,a^\vee\rangle=m$ when
$a=2+k\delta$ for some $k\in {\mathbb{Z}}$.
The lemma is now immediate when $m\in
{\mathbb{Z}}_+$. For $m\in -{\mathbb{N}}$, we first observe
that the coefficient of $x^{-m}$ in the expansion of
${\mathcal{R}}(a)(x^m)$ in terms of monomials is zero. Indeed, the
coefficient of $x^{-m}$ in
\[t_a^{-1}(1-t_at_{a/2}x^{a/2})
(1+t_at_{a/2}^{-1}x^{a/2})D_a(x^m)\]
is $-t_aq^{-mk}$, which cancels with the coefficient of
$t_as_a(x^m)=t_aq^{-mk}x^{-m}$.
Hence the highest order term of ${\mathcal{R}}(a)(x^m)$
is $t_a^{-1}x^m$ when $m\in -{\mathbb{N}}$. This completes
the proof of the lemma.
\end{proof}

\subsection {}

%3.4
Lemma 3.3 implies the following
triangularity property of the Cherednik-Dunkl operator $Y$.
Set $\gamma_m:=k_0^{\epsilon(m)}k_1^{\epsilon(m)}q^m$ for
$m\in {\mathbb{Z}}$.
\begin{prop}
For all $m\in {\mathbb{Z}}$, we have
\[ Y(x^m)=\gamma_mx^m+\, \hbox{\textup{lower order terms w.r.t.}}\,
\preceq.
\]
\end{prop}
\begin{proof}
Observe that
$Y=T_1T_0={\mathcal{R}}(a_1)s_1{\mathcal{R}}(a_0)s_0
={\mathcal{R}}(a_1){\mathcal{R}}(s_1(a_0))\tau(1)$.
Now $s_1(a_0)=2+\delta$ and $\tau(1)(x^m)=q^mx^m$ (see 2.10), so
the
proposition follows from lemma 3.3.
\end{proof}

\subsection {}

%3.5
The diagonal terms $\gamma_m$ ($m\in {\mathbb{Z}}$) of the triangular
operator $Y$ are pair-wise different by the
generic conditions 2.11
on $q$ and on the multiplicity function $\underline{t}$.
Hence proposition 3.4 leads immediately to the
following proposition (compare with Sahi \cite[section 6]{S} for the
higher rank setting).
\begin{prop}
There exists a unique basis
$\{P_m(\cdot)=P_m(\cdot;\underline{t};q)\,\, | \,\, m\in {\mathbb{Z}}\}$ of
${\mathcal{A}}$ such that
\begin{enumerate}
\item{}
$P_m(x)=x^m+$ lower order terms with respect to $\preceq$,
\item{}
$Y(P_m)=\gamma_mP_m$
\end{enumerate}
for all $m\in {\mathbb{Z}}$.
\end{prop}

\begin{Def}
The Laurent polynomial $P_m=P_m(\cdot;\underline{t};q)$
\textup{(}$m\in {\mathbb{Z}}$\textup{)}
is called the monic, non-sym\-me\-tric Askey-Wilson polynomial of degree
$m$.
\end{Def}
We will justify this terminology in section 5, where we
relate the non-symmetric Askey-Wilson polynomials with the well-known
symmetric Askey-Wilson polynomials by a kind of symmetrization procedure.

%%%%%%%%%%%%%%%%%%%%%%%%%%%%%%%%%%%%%%
%%                                  %%
%%                                  %%
%%  The fundamental representation  %%
%%                                  %%
%%                                  %%
%%%%%%%%%%%%%%%%%%%%%%%%%%%%%%%%%%%%%%

\section{The fundamental representation}

\subsection {}

%4.1
In the previous section we have diagonalized the action of the
``translation part'' ${\mathbb{C}}[Y^{\pm 1}]$ of the affine Hecke
algebra $H=H(R;k_0,k_1)$ under the fundamental
representation $\pi_{\underline{t},q}$ (see 2.19). The corresponding
eigenfunctions are exactly the non-symmetric Askey-Wilson
polynomials.
Since $H$ is generated as algebra by $Y$ and the
difference-reflection operator $T_1$, see 2.15, it suffices to understand
the action of $T_1$ on the non-symmetric Askey-Wilson polynomials in
order to completely decompose ${\mathcal{A}}$ as an $H$-module.
Recall the notation
$\gamma_m=k_0^{\epsilon(m)}k_1^{\epsilon(m)}q^m$
($m\in {\mathbb{Z}}$) for the
eigenvalues of $Y$.

\begin{prop}
For $m\in {\mathbb{Z}}$ we have
\[ T_1P_m=\alpha_mP_m+\beta_mP_{-m},
\]
with
\[ \alpha_m=\frac{(k_1^{-1}-k_1)\gamma_m^2
+(k_0^{-1}-k_0)\gamma_m}{1-\gamma_m^2}.
\]
If $m\in -{\mathbb{N}}$ then $\beta_m=k_1$, and
\[
\beta_m=k_1\prod_{\xi=\pm 1}
\frac{(1+k_0k_1^{-1}\gamma_m^{\xi})(1-k_0^{-1}k_1^{-1}\gamma_m^{\xi})}
{(1-\gamma_m^{2\xi})}
\]
if $m\in {\mathbb{Z}}_+$.
\end{prop}
\begin{proof}
The formula for $m=0$ reduces to $T_1(P_0)=k_1P_0$, which is clear.
For $0\not=m\in {\mathbb{Z}}$,
we derive from Lusztig's formula 2.16 and from the definition 3.5 of the
non-symmetric Askey-Wilson polynomials that for all
$f(Y)\in {\mathbb{C}}[Y^{\pm 1}]$,
\[ \bigl(f(Y)-f(\gamma_{-m})\bigr)T_1P_m=
\alpha_m\bigl(f\bigl(\gamma_m\bigr)-
f\bigl(\gamma_{-m}\bigr)\bigr)P_m
\]
with $\alpha_m$ as given in the statement of the proposition.
Since $\gamma_m$ ($m\in {\mathbb{Z}}$) are mutually
different by the conditions 2.11 on the parameters,
we derive from proposition 3.5 that $T_1P_m=\alpha_mP_m+\beta_mP_{-m}$ for
some
$\beta_m$.

If $m\in -{\mathbb{N}}$, then we have $x^m\prec s_1(x^m)=x^{-m}$. Combined
with
the formula
$T_1=s_1{\mathcal{R}}(a_1)^{-1}+k_1-k_1^{-1}$ (see 3.1) and with the
triangularity of ${\mathcal{R}}(a_1)$ (see lemma 3.3),
we obtain that the coefficient of $x^{-m}$ in
the expansion of $T_1(x^m)$ with respect to the basis of monomials is equal
to $k_1$. By the definition 3.5 of the non-symmetric Askey-Wilson
polynomials, we conclude that $\beta_m=k_1$ for $m\in -{\mathbb{N}}$.

Now act by $T_1$ on both sides of the formula $T_1P_m=\alpha_mP_m+
\beta_mP_{-m}$ and use the quadratic relation for $T_1$,
see 2.12. It follows that the $\alpha_m$ and the $\beta_m$ satisfy
the
relation
\[ \beta_m\beta_{-m}=(k_1-\alpha_m)(k_1^{-1}+\alpha_m),\qquad
0\not=m\in {\mathbb{Z}}.
\]
This allows us to compute $\beta_m$ with $m\in {\mathbb{N}}$ from the known
expressions for $\alpha_m$ and $\beta_{-m}$, which yields the
desired result.
\end{proof}

\subsection {}

%4.2
A uniform formula for the action of $T_1$ on the non-symmetric
Askey-Wilson polynomials can be obtained by renormalizing the
non-symmetric Askey-Wilson polynomials in a suitable way. A
natural renormalization, together with a new proof of proposition 4.1,
is given in section 10.

\subsection {}

%4.3
As a consequence of proposition 4.1, we can compute the action of
the non-affine intertwiner $S_1$ (see 2.17)
on the non-symmetric Askey-Wilson polynomials explicitly.
\begin{cor}
We have $S_1(P_m)=
\bigl(\gamma_m-\gamma_{-m}\bigr)\beta_mP_{-m}$ for $m\in
{\mathbb{Z}}$,
where $\beta_m$ is as in proposition 4.1.
\end{cor}
\begin{proof}
By proposition 4.1
and by the definition 3.5 of the non-symmetric Askey-Wilson
polynomial, we have
\begin{equation*}
\begin{split}
S_1P_m=&(T_1Y-YT_1)P_m=\bigl(\gamma_m-Y\bigr)T_1P_m\\
=&\bigl(\gamma_m-Y\bigr)\bigl(\alpha_mP_m+\beta_mP_{-m}\bigr)
=\bigl(\gamma_m-\gamma_{-m}\bigr)\beta_mP_{-m}.
\end{split}
\end{equation*}
\end{proof}

\subsection {}

%4.4
We set ${\mathcal{A}}(0)=\hbox{span}\{P_0\}$ and
${\mathcal{A}}(m)=\hbox{span}\{P_m, P_{-m}\}$ for $m\in
{\mathbb{N}}$.

\begin{thm}
{\bf (i)}
The representation $(\pi_{\underline{t},q},H(R;k_0,k_1))$ is faithful.

{\bf (ii)}
The center ${\mathcal{Z}}(H)$ of $H=H(R;k_0,k_1)$ is equal to
${\mathbb{C}}[Y^{\pm 1}]^W={\mathbb{C}}[Y+Y^{-1}]$.

{\bf (iii)}
The decomposition ${\mathcal{A}}=
\oplus_{m\in {\mathbb{Z}}_+}{\mathcal{A}}(m)$ is the
multiplicity-free, irreducible decomposition of ${\mathcal{A}}$
as $(\pi_{\underline{t},q},H)$-module. It is also the
decomposition of ${\mathcal{A}}$ in isotypical components for
the action of the center, where the central character of
${\mathcal{A}}(m)$ is given by
$\chi_m(f(Y))=f\bigl(\gamma_m\bigr)$ for
$f(Y)\in {\mathcal{Z}}(H)={\mathbb{C}}[Y^{\pm 1}]^W$.
\end{thm}
\begin{proof}

{\bf (i)}
Suppose that $h=f(Y)+T_1g(Y)$ acts as zero on ${\mathcal{A}}$ under
the representation $\pi_{\underline{t},q}$,
where $f,g\in {\mathcal{A}}$. Let $h$
act on the non-symmetric Askey-Wilson polynomials $P_m$ ($m<0$)
and use proposition 4.1 together with the fact that the coefficients
$\beta_m$ ($m<0$) in 4.1 are non-zero. Then we conclude that
$g(\gamma_m)=0$ for all $m\in -{\mathbb{N}}$. By the
conditions 2.11 on the parameters, this implies $g=0$ in ${\mathcal{A}}$.
But $h=f(Y)$ acting on $P_m$ shows that
$f\bigl(\gamma_m\bigr)=0$ for all $m$,
hence $f=0$ in ${\mathcal{A}}$. Combined with 2.15, this shows that
$\pi_{\underline{t},q}$ is faithful.

{\bf (ii)}
Clearly any element from ${\mathbb{C}}[Y+Y^{-1}]$ commutes with
${\mathbb{C}}[Y^{\pm 1}]$, but also with $T_1$ by Lusztig's formula
2.16. Hence 2.15 gives ${\mathbb{C}}[Y+Y^{-1}]\subset
{\mathcal{Z}}(H)$.
Suppose $0\not=h=f(Y)+T_1g(Y)\in {\mathcal{Z}}(H)$,
where $f, g\in {\mathcal{A}}$.
Then $h$ acts as a constant on each of the
$P_m$ ($m\in {\mathbb{Z}}$). In view of proposition 4.1, this implies
that $g(\gamma_m)=0$ for all $m\not=0$, hence $g=0$ in ${\mathcal{A}}$.
By corollary 4.3 we then have for $m\not=0$,
\[
f\bigl(\gamma_m\bigr)S_1P_m=S_1(hP_m)=h(S_1P_m)
=f\bigl(\gamma_m^{-1}\bigr)S_1P_m.
\]
Furthermore, $S_1P_m\not=0$ by the conditions 2.11 on the
parameters. Hence $f(\gamma_m)=f(\gamma_m^{-1})$ for $0\not=m\in
{\mathbb{Z}}$,
i.e. $h=f(Y)\in {\mathbb{C}}[Y+Y^{-1}]$.

{\bf (iii)} The second statement follows directly from
proposition 3.5 and from the fact that the central character values
$\chi_m(Y+Y^{-1})=\gamma_m+\gamma_m^{-1}$ ($m\in
{\mathbb{Z}}_+$) are pair-wise different by the conditions 2.11 on the
parameters.
For the first statement, it then suffices to show that
${\mathcal{A}}(m)$ ($m\in {\mathbb{Z}}_+$) are irreducible
$H$-modules. This follows without difficulty
from 2.15, proposition 4.1, corollary 4.3 and the fact
that $\beta_m\not=0$ for all $0\not=m\in {\mathbb{Z}}$ by the
conditions 2.11 on the parameters.
\end{proof}

%%%%%%%%%%%%%%%%%%%%%%%%%%%%%%%%%%%%%%%%%%%%%%%%%%%%%%%%%%%%%%%%%%%%
%%                                                                %%
%%    The (anti-)symmetric Askey-Wilson polynomials               %%
%%                                                                %%
%%%%%%%%%%%%%%%%%%%%%%%%%%%%%%%%%%%%%%%%%%%%%%%%%%%%%%%%%%%%%%%%%%%%

\section{The (anti-)symmetric Askey-Wilson polynomials}

\subsection {}

%5.1
In the present rank one setting, the representation theory of the
underlying two-dimensional Hecke algebra $H_0=H_0(k_1)$ is
extremely
simple: the trivial representation $\chi_+$ and the alternating
representation $\chi_-$ exhaust its irreducible representations,
where
$\chi_\pm$ are uniquely determined by $\chi_{\pm}(T_1)=\pm k_1^{\pm 1}$.
The corresponding mutually orthogonal, primitive idempotents are
given
by
\[ C_+=\frac{1}{1+k_1^2}\bigl(1+k_1T_1\bigr),\qquad
C_-=\frac{1}{1+k_1^{-2}}\bigl(1-k_1^{-1}T_1\bigr).
\]
So $\{C_-,C_+\}$ is a partition of the unity for $H_0$.
In particular, we have $C_-+C_+=1$.

\subsection {}

%5.2
The partition of the unity of $H_0$ introduced in 5.1 gives the
decomposition ${\mathcal{A}}={\mathcal{A}}_-\oplus
{\mathcal{A}}_+$ of ${\mathcal{A}}$
in isotypical components for the action of
$(\pi_{\underline{t},q}|_{H_0}, H_0)$, where
${\mathcal{A}}_{\pm}=C_{\pm}{\mathcal{A}}$.
Observe that ${\mathcal{A}}_{\pm}$ consists of the Laurent
polynomials $f\in {\mathcal{A}}$ which satisfy
$(T_1\mp k_1^{\pm 1})f=0$. By the explicit expression 2.19 for the
difference-reflection operator $T_1$
we have
$T_1-k_1=\phi_1(x)(s_1-\hbox{id})$ for some non-zero rational
function $\phi_1(x)$, so that ${\mathcal{A}}_+$ coincides with the algebra
${\mathcal{A}}^W={\mathbb{C}}[x+x^{-1}]$ consisting of $W$-invariant Laurent
polynomials. We call ${\mathcal{A}}_-$ the subspace of
{\it anti-symmetric} Laurent polynomials.
The decomposition $f=C_-f+C_+f$ for $f\in {\mathcal{A}}$ is its
unique
decomposition as a sum of an anti-symmetric and a symmetric Laurent
polynomial.

\subsection {}

%5.3
The irreducible $H$-module ${\mathcal{A}}(m)\subset
{\mathcal{A}}$
decomposes under the action of $H_0$ by
${\mathcal{A}}(m)={\mathcal{A}}_-(m)\oplus {\mathcal{A}}_+(m)$,
where
${\mathcal{A}}_{\pm}(m)=C_{\pm}{\mathcal{A}}(m)$.

\begin{prop}
{\bf (i)} Let $m\in {\mathbb{Z}}_+$, then
$\hbox{dim}({\mathcal{A}}_{+}(m))=1$.
More precisely, there exists a unique
$P_m^+\in {\mathcal{A}}_+(m)$ of the form $P_m^+(x)=x^m+$
lower order terms with respect to $\preceq$. In terms of
non-symmetric
Askey-Wilson polynomials, we have
\[P_m^+=P_m+\frac{(1+k_0k_1^{-1}\gamma_m)(1-k_0^{-1}k_1^{-1}\gamma_m)}
{(1-\gamma_m^2)}P_{-m},\qquad m\in {\mathbb{Z}}_+.
\]

{\bf (ii)} We have ${\mathcal{A}}_-(0)=\{0\}$ and
$\hbox{dim}({\mathcal{A}}_-(m))=1$ for $m\in {\mathbb{N}}$. More precisely,
there exists for all $m\in {\mathbb{N}}$ a unique $P_m^-\in
{\mathcal{A}}_-(m)$ of the form $P_m^-(x)=x^m+$ lower order terms
with respect to $\preceq$.
In terms of non-symmetric Askey-Wilson polynomials, we have
\[P_m^-=P_m-
\frac{(1+k_0k_1^{-1}\gamma_m^{-1})(1-k_0^{-1}k_1^{-1}\gamma_m^{-1})}
{(1-\gamma_m^{-2})}P_{-m},
\qquad m\in {\mathbb{N}}.
\]
\end{prop}
\begin{proof}
The statements for $m=0$ are immediate since $T_1(1)=k_1 1$, where
$1\in {\mathcal{A}}$ is the Laurent polynomial identically equal to one.
For $m\in {\mathbb{N}}$,
we can write $C_{\pm}P_m\in {\mathcal{A}}_{\pm}(m)$ explicitly as
\begin{equation*}
\begin{split}
C_{\pm}P_m&=\frac{1\pm k_1^{\pm 1}\alpha_m}{1+k_1^{\pm 2}}P_m\pm
\frac{k_1^{\pm 1}\beta_m}{1+k_1^{\pm 2}}P_{-m}\\
&=
k_1^{\pm 1}\frac{(k_1^{\mp 1}\pm \alpha_m)}{(1+k_1^{\pm 2})}
\bigl(P_m\pm k_1^{-1}(k_1^{\pm 1}\mp \alpha_m)P_{-m}\bigr)
\end{split}
\end{equation*}
in view of (the proof of) proposition 4.1.
Observe that the coefficient of $P_m$ is non-zero by the
conditions 2.11 on the parameters. In particular,
${\mathcal{A}}_{\pm}(m)=\hbox{span}\{C_{\pm}P_m\}$
are one-dimensional subspaces for all
$m\in {\mathbb{N}}$. Dividing out the non-zero coefficient
of $P_m$ in the expansion of $C_{\pm}P_m$ and using proposition 3.5,
we conclude that
there exist unique elements $P_m^{\pm}\in {\mathcal{A}}_{\pm}(m)$ ($m\in
{\mathbb{N}}$) satisfying $P_m^{\pm}(x)=x^m+$ lower order terms w.r.t.
$\preceq$.
The explicit formulas for $P_m^{\pm}$ in terms of non-symmetric
Askey-Wilson polynomials follow now by substituting the explicit
expression for $\alpha_m$ in the above
expansion of $C_{\pm}P_m$, see proposition 4.1.
\end{proof}

\begin{Def}

{\bf (i)} The polynomial $P_m^+=P_m^+(\cdot;\underline{t};q)
\in {\mathcal{A}}_+$ \textup{(}$m\in
{\mathbb{Z}}_+$\textup{)} is called the monic, symmetric Askey-Wilson
polynomial of degree $m$.

{\bf (ii)} The polynomial $P_m^-=P_m^-(\cdot;\underline{t};q)
\in {\mathcal{A}}_-$ \textup{(}$m\in {\mathbb{N}}$\textup{)}
is called the monic, anti-symmetric Askey-Wilson polynomial of degree
$m$.
\end{Def}

Askey and Wilson \cite{AW} defined a very general
family of basic hypergeometric orthogonal polynomials which are
nowadays known as the Askey-Wilson polynomials. In 5.9 we justify our
terminology for the Laurent polynomials $P_m^{\pm}$
by showing that the $P_m^+$ ($m\in {\mathbb{Z}}$) coincide with
the Askey-Wilson polynomials as
defined in \cite{AW}.

\subsection {}

%5.4
We can also express the non-symmetric Askey-Wilson polynomial
in terms of the symmetric Askey-Wilson polynomial in the following
way.
\begin{lem}
We have $P_0=P_0^+$, and for $m\in {\mathbb{N}}$,
\begin{equation*}
\begin{split}
P_m&=\frac{1}{\gamma_m-\gamma_{-m}}\bigl(Y-\gamma_{-m}\bigr)P_m^+,\\
P_{-m}&=\frac{\gamma_m}{(1+k_0k_1^{-1}\gamma_m)(1-k_0^{-1}k_1^{-1}\gamma_m)}
\bigl(Y-\gamma_m\bigr)P_m^+.
\end{split}
\end{equation*}
\end{lem}
\begin{proof}
The statement for $m=0$ is trivial. For $m\in
{\mathbb{N}}$ the formulas follow directly from proposition 3.5 and from
the expansion of $P_m^+$ as linear combination of non-symmetric
Askey-Wilson polynomials, see proposition 5.3.
\end{proof}

\subsection {}

%5.5
Observe that the affine Weyl group
${\mathcal{W}}$ acts on ${\mathcal{A}}$ by algebra
automorphisms (see 2.10).
This action can be uniquely extended to an action (by automorphisms)
of $\mathcal{W}$ on the rational functions ${\mathbb{C}}(x)$ in the
indeterminate $x$.
Since $|q|\not=1$ (see 2.11), we have
\[ \bigoplus_{w\in {\mathcal{W}}}{\mathbb{C}}(x)w=\bigoplus_{m\in
{\mathbb{Z}},
\sigma\in W}{\mathbb{C}}(x)\tau(m)\sigma
\]
as a subalgebra of
$\hbox{End}_{\mathbb{C}}\bigl({\mathbb{C}}(x)\bigr)$.

\subsection {}

%5.6
Any $X\in \hbox{Im}(\pi_{\underline{t},q})\subset
\hbox{End}_{\mathbb{C}}\bigl({\mathcal{A}})$ can be uniquely
written as a finite ${\mathbb{C}}(x)$-linear
combination of the automorphisms $w\in {\mathcal{W}}$
of ${\mathcal{A}}$. By 5.5, we may regard $X$ as a linear endomorphism of
${\mathbb{C}}(x)$. We thus have a unique decomposition
$X=X_-s_1+X_+$ where $X_{\pm}\in \bigoplus_{m\in
{\mathbb{Z}}}{\mathbb{C}}(x)\tau(m)$ are $q$-difference operators
with rational coefficients. We write $X_{sym}:=X_-+X_+$, so that
$Xf=X_{sym}f$ for all $f\in {\mathcal{A}}_+={\mathcal{A}}^W$.

\subsection {}

%5.7
In order to make the connection between the symmetric Askey-Wilson
polynomials
$P_m^+$ ($m\in {\mathbb{Z}}$) and the four parameter family of
Askey-Wilson polynomials as originally defined in \cite{AW},
it is convenient to reparametrize the multiplicity function
$\underline{t}\simeq
(u_0,u_1,k_0,k_1)$ (see 2.11) by
\[ a=k_1u_1,\quad b=-k_1u_1^{-1},\quad c=q^{\frac{1}{2}}k_0u_0,\quad
d=-q^{\frac{1}{2}}k_0u_0^{-1}.
\]

\subsection {}

%5.8
Using the parameters $a,b,c,d$ (see 5.7),
we can give the following explicit expression
for the $q$-difference operator $(Y+Y^{-1})_{sym}$, see \cite{N}
for the higher rank result.
\begin{prop}
We have
\[ \bigl(Y+Y^{-1}\bigr)_{sym}=A(x)\bigl(\tau(1)-1\bigr)+A(x^{-1})
\bigl(\tau(-1)-1\bigr)+
k_0k_1+k_0^{-1}k_1^{-1},
\]
with
\[ A(x)=k_0^{-1}k_1^{-1}\frac{(1-ax)(1-bx)(1-cx)(1-dx)}{(1-x^2)(1-qx^2)}.
\]
\end{prop}
\begin{proof}
We write $T_i=k_i+\phi_i(x)(s_i-1)$ with
\[
\phi_0(x)=k_0^{-1}\frac{(1-cx^{-1})(1-dx^{-1})}{(1-qx^{-2})},\quad
\phi_1(x)=k_1^{-1}\frac{(1-ax)(1-bx)}{(1-x^2)}
\]
for the difference-reflection operators associated with $S$, see 2.19.
Since $Y=T_1T_0$, $s_0=t(-1)s_1$ and $T_i^{-1}=T_i+k_i^{-1}-k_i$ for
$i=0,1$, we have
\[
(Y+Y^{-1})_{sym}=B(x)(\tau(1)-1)+C(x)(\tau(-1)-1)+D(x)
\]
for unique coefficients $B,C,D\in {\mathbb{C}}(x)$. Observe that
$D(x)=(Y+Y^{-1})(1)=k_0k_1+k_0^{-1}k_1^{-1}$ where $P_0=1\in
A_+$ is the Laurent polynomial identically equal to one, since $T_i(1)=k_i
1$
for $i=0,1$. To compute the coefficient $B(x)$ (respectively $C(x)$), we
need to
compute the coefficient of $\tau(1)$ (respectively $\tau(-1)$) in
$(Y+Y^{-1})_{sym}$.
Recall that $s_0=s_1\tau(1)=\tau(-1)s_1$, so that the
$\tau(1)$-term (respectively $\tau(-1)$-term) of $Y_{sym}$ has coefficient
$\phi_1(x)\phi_0(x^{-1})=A(x)$ (respectively has coefficient
$(k_1-\phi_1(x))\phi_0(x)$). The $\tau(1)$-term
(respectively $\tau(-1)$-term) of $(Y^{-1})_{sym}$ is zero
(respectively has coefficient $\phi_0(x)\phi_1(qx^{-1})+
\phi_0(x)(k_1^{-1}-\phi_1(qx^{-1}))=
k_1^{-1}\phi_0(x)$). Adding the contributions, we see that $B(x)=A(x)$
and that
$C(x)=\phi_0(x)\bigl(k_1+k_1^{-1}-\phi_1(x)\bigr)=A(x^{-1})$,
which completes the proof of the proposition.
\end{proof}
The second order $q$-difference operator $L=(Y+Y^{-1})_{sym}$ is
called the Askey-Wilson second-order $q$-difference operator, cf.
\cite[(5.7)]{AW}.

\subsection {}

%5.9
The symmetric Askey-Wilson polynomial $P_m^+$ ($m\in {\mathbb{Z}}_+$)
lies in the irreducible $H(R;k_0,k_1)$-module
${\mathcal{A}}(m)$,
hence the central element $Y+Y^{-1}\in {\mathcal{Z}}(H)$ acts on
$P_m^+$
as the scalar $\gamma_m+\gamma_{m}^{-1}$, see theorem 4.4. Combined
with  proposition 5.8, we conclude
that $P_m^+$ is an eigenfunction of the Askey-Wilson
second-order
$q$-difference operator $L$ with eigenvalue $\gamma_m+\gamma_m^{-1}$.
The eigenvalues $\gamma_m+\gamma_m^{-1}$ ($m\in {\mathbb{Z}}_+$)
are mutually different by the conditions 2.11 on the parameters and
$\{P_m^+\}_{m\in {\mathbb{Z}}_+}$ is a linear basis of
${\mathcal{A}}_+={\mathcal{A}}^W$, so that $P_m^+$
is the unique $W$-invariant Laurent polynomial satisfying
$LP_m^+=(\gamma_m+\gamma_m^{-1})P_m^+$.
A comparison with \cite[(5.7)]{AW} yields the following result.
\begin{thm}
The $W$-invariant Laurent polynomial $P_m^+$ \textup{(}$m\in
{\mathbb{Z}}_+$\textup{)}
coincides with the monic Askey-Wilson polynomial
of degree $m$ as described in \cite{AW}. In particular,
we have in terms of basic hypergeometric series,
\[
P_m^+(x)=\frac{\bigl(ab,ac,ad;q\bigr)_m}{a^m\bigl(abcdq^{m-1};q\bigr)_m}
{}_4\phi_3\left(\begin{matrix} &q^{-m},  q^{m-1}abcd, ax, ax^{-1} \\
& ab, ac, ad \end{matrix}; q,q\right).
\]
\end{thm}

\subsection {}

%5.10
Theorem 5.9 and lemma 5.4 can be used to write the
non-symmetric and the anti-symmetric Askey-Wilson polynomials as
a sum of two terminating balanced ${}_4\phi_3$'s.
It is convenient to write $P_m^+(x)=P_m^+(x;a,b,c,d)$ for the
symmetric Askey-Wilson polynomial $P_m^+(x)=P_m^+(x;\underline{t};q)$
when we want to emphasize the
dependence of $P_m^+$
on the (reparametrized) multiplicity function $(a,b,c,d)$, see 5.7.

\begin{prop}

{\bf (i)} For $m\in {\mathbb{Z}}_+$ we have
\begin{equation*}
\begin{split}
&P_m(x)=q^m\frac{(1-abcdq^{m-1})}{(1-abcdq^{2m-1})}P_m^+(x;a,b,c,d)\\
&\,+q^{(m-1)/2}\frac{(1-cx^{-1})(1-dx^{-1})x(1-q^m)}{(1-abcdq^{2m-1})}
P_{m-1}^+(q^{-1/2}x;q^{1/2}a,q^{1/2}b,q^{1/2}c,q^{1/2}d),
\end{split}
\end{equation*}
where the second term should be read as zero when $m=0$.

{\bf (ii)} For $m\in {\mathbb{N}}$ we have
\begin{equation*}
\begin{split}
P_{-m}(x)&=\frac{1}{(1-cdq^{m-1})}P_m^+(x;a,b,c,d)\\
&-q^{(m-1)/2}\frac{(1-cx^{-1})(1-dx^{-1})x}{(1-cdq^{m-1})}
P_{m-1}^+(q^{-1/2}x;q^{1/2}a,
q^{1/2}b,q^{1/2}c,q^{1/2}d).
\end{split}
\end{equation*}

{\bf (iii)} For $m\in {\mathbb{N}}$ we have
\begin{equation*}
\begin{split}
&P_m^-(x)=\frac{(1-abcdq^{m-1})}{ab(1-cdq^{m-1})}
P_m^+(x;a,b,c,d)\\
&\,+q^{(m-1)/2}\frac{(1-cx^{-1})(1-dx^{-1})x(ab-1)}{ab(1-cdq^{m-1})}
P_{m-1}^+(q^{-1/2}x;q^{1/2}a,q^{1/2}b,q^{1/2}c,q^{1/2}d).
\end{split}
\end{equation*}
\end{prop}
\begin{proof}
Recall the rational function $\phi_0(x)$ defined in the proof of
proposition 5.8.
The proof of proposition 5.8 shows that
$(Y^{-1})_{sym}=k_1^{-1}\phi_0(x)\bigl(\tau(-1)-1\bigr)+k_0^{-1}k_1^{-1}$.

{\bf (i)} The formula for $m=0$ is trivial. Let $m\in {\mathbb{N}}$.
By lemma 5.4, we have
$P_m=(\gamma_m-\gamma_{-m})^{-1}(\gamma_m-(Y^{-1})_{sym})P_m^+$.
In view of the explicit formula for the $q$-difference operator
$(Y^{-1})_{sym}$, we need to write $(t(-1)-1)P_m^+$ as a terminating
balanced
${}_4\phi_3$. By a direct computation using the explicit expression of
$P_m^+$ in terms of a terminating balanced ${}_4\phi_3$, see theorem 5.9,
we have
\begin{equation*}
\begin{split}
\bigl((\tau(-1)-1)P_m^+(.;a,b,c,d)\bigr)(x)=
&(q^{1/2}x^{-1}-q^{-1/2}x)(q^{m/2}-q^{-m/2})\\
&\quad.P_{m-1}^+(q^{-1/2}x;q^{1/2}a, q^{1/2}b,q^{1/2}c,q^{1/2}d),
\end{split}
\end{equation*}
cf. \cite[(5.6)]{AW} or \cite[(7.7.7)]{GR}. This leads to the desired
result.

{\bf (ii)} We have
$P_{-m}(x)=(1+k_0k_1^{-1}\gamma_m)^{-1}(1-k_0^{-1}k_1^{-1}\gamma_m)^{-1}
(1-\gamma_m(Y^{-1})_{sym})P_m^+$ by lemma 5.4. The proof is now similar
to the proof of {\bf (i)}.

{\bf (iii)} This follows from {\bf (i)} and {\bf (ii)}, together
with proposition 5.3.
\end{proof}

%%%%%%%%%%%%%%%%%%%%%%%%%%%%%%%%%%%%%%%%%%%%%%%%%%%%%%%%%%%%%%%%%%%%
%%                                                                %%
%%                 (Bi-)orthogonality relations                   %%
%%                                                                %%
%%                                                                %%
%%%%%%%%%%%%%%%%%%%%%%%%%%%%%%%%%%%%%%%%%%%%%%%%%%%%%%%%%%%%%%%%%%%%

\section{(Bi-)orthogonality relations}

\subsection {}
%6.1
We assume in this section that the multiplicity function $\underline{t}$
and the deformation parameter $q$ satisfy the additional conditions
that $0<|q|<1$ and that $ef\not\in q^{{\mathbb{Z}}}$ for all
$e,f\in \{a,b,c,d\}$.

\subsection {}

%6.2
Let $C\subset {\mathbb{C}}$ be a continuous rectifiable Jordan curve
such that $aq^k, bq^k,cq^k, dq^k$ ($k\in {\mathbb{Z}}_+$) are in the
interior of $C$ and such that
$a^{-1}q^{-k}, b^{-1}q^{-k}, c^{-1}q^{-k}, d^{-1}q^{-k}$
($k\in {\mathbb{Z}}_+$) are in the exterior of $C$.
By the conditions 6.1 on the parameters, such a contour
exists. We give $C$ the counterclockwise orientation.
Let $\langle .,. \rangle=
\langle .,. \rangle_{\underline{t},q}$ and
$\bigl( .,. \bigr)=\bigl( .,. \bigr)_{\underline{t},q}$ be the
bilinear forms on ${\mathcal{A}}$ defined by
\[\langle f,g\rangle =
\frac{1}{2\pi i}\int_{x\in C}f(x)g(x^{-1})\Delta(x)\frac{dx}{x},
\qquad
\bigl(f,g\bigr) =\frac{1}{2\pi i}\int_{x \in C}
f(x)g(x^{-1})\Delta_+(x)\frac{dx}{x},
\]
where the weight functions $\Delta(x)=\Delta(x;\underline{t};q)$
and $\Delta_+(x)=
\Delta_+(x;\underline{t};q)$ are given by the infinite products
\begin{equation*}
\begin{split}
\Delta(x)&=\prod_{a\in R^+}
\frac{\bigl(1-x^a\bigr)}{\bigl(1-t_at_{a/2}x^{a/2}\bigr)
\bigl(1+t_at_{a/2}^{-1}x^{a/2}\bigr)},\\
\Delta_+(x)&=\prod_{a\in R: a(0)\geq 0}\frac{\bigl(1-x^a\bigr)}
{\bigl(1-t_at_{a/2}x^{a/2}\bigr)
\bigl(1+t_at_{a/2}^{-1}x^{a/2}\bigr)}.
\end{split}
\end{equation*}
The conditions 6.1 on the parameters
ensure that the weight functions are well-defined.
In terms of $q$-shifted factorials, we can rewrite the
weight function $\Delta_+(x)$ as
\[ \Delta_+(x)=\frac{\bigl(x^2,x^{-2};q\bigr)_{\infty}}
{\bigl(ax,ax^{-1},bx,bx^{-1},cx,cx^{-1},dx,dx^{-1};q\bigr)_{\infty}}
\]
using 5.7, 2.7 and 2.10. Hence $\Delta_+(\cdot)$ coincides with the
weight function of the orthogonality measure of
the symmetric Askey-Wilson polynomials as defined in \cite{AW}
(see also proposition 6.9).
Observe that
$\Delta(x)=\alpha(x)\Delta_+(x)$ with $\alpha(x)=\alpha(x;k_1,u_1)$
given by
\[\alpha(x)=\frac{(1-k_1u_1x^{-1})(1+k_1u_1^{-1}x^{-1})}{(1-x^{-2})}=
\frac{(1-ax^{-1})(1-bx^{-1})}{(1-x^{-2})}.
\]

\subsection {}

%6.3
Using Cauchy's theorem we can rewrite
$\langle .,. \rangle$ and $\bigl( .,. \bigr)$ as an integral
over the unit circle $T$ in the complex plane plus a finite sum
of residues of the integrand.
The residues of the weight functions $\Delta(\cdot)$
and $\Delta_+(\cdot)$ can be computed explicitly, see
\cite[section 2]{AW} or \cite[section 7.5]{GR} for more details.

\subsection {}

%6.4
Observe that the factor $\alpha(x)$ in the weight function
$\Delta(x)$ satisfies the identity
$\alpha(x)+\alpha(x^{-1})=1-ab$.
Since $\Delta_+(x)$ is furthermore invariant under $x\mapsto x^{-1}$,
we see that the restrictions of the bilinear forms
$\langle .,. \rangle$ and $\bigl( .,. \bigr)$ to
${\mathcal{A}}^W$ coincide up to a constant:
\begin{lem}
For $f,g\in {\mathcal{A}}^W$, we have
$\langle
f,g\rangle=\frac{1}{2}(1-ab)\,\bigl(f,g\bigr)=\frac{1}{2}(1+k_1^2)\,
\bigl(f,g\bigr)$.
\end{lem}

\subsection {}

%6.5
Let $T$ be a linear endomorphism of ${\mathcal{A}}$. Then there
exists at most one linear endomorphism $T^*$ of ${\mathcal{A}}$
such that $\langle Tf,g\rangle =\langle f,T^*g\rangle$ for all $f,g\in
{\mathcal{A}}$, since the bilinear form $\langle .,. \rangle$
is non-degenerate. If $T^*$ exists, then
we call $T^*$ the {\it adjoint} of $T$ with respect
to $\langle .,. \rangle$.

\subsection {}

%6.6
We write $T_0^\prime, T_1^\prime$ for the difference-reflection
operators associated to $S$ with respect to inverse parameters
$(\underline{t}^{-1},q^{-1})$, where $\underline{t}^{-1}$ is the
multiplicity function $(t_a^{-1})_{a\in S}$.
More precisely, $T_i^\prime$ is the
image of the fundamental generator $T_i\in H(R;k_0^{-1},k_1^{-1})$
under the (faithful) representation $\pi_{\underline{t}^{-1},q^{-1}}$,
see 2.19.
Furthermore, we set $Y^\prime=T_1^{\prime}T_0^\prime$
for the associated Cherednik-Dunkl operator.
\begin{prop}
The adjoint of the difference-reflection operator $T_i$
\textup{(}$i=0,1$\textup{)}
and of the Dunkl operator $Y$ exists.
More precisely, we have $T_i^*=(T_i^\prime)^{-1}$
\textup{(}$i=0,1$\textup{)} and
$Y^*=(Y^\prime)^{-1}$.
\end{prop}
\begin{proof}
We use the notation $T_i=k_i+\phi_i(x)(s_i-1)$ and
$T_i^\prime=k_i^{-1}+\phi_i^\prime(x)(s_i^\prime-1)$ ($i=0,1$)  for the
difference reflection operator with respect to the
parameters $(\underline{t},q)$ and $(\underline{t}^{-1},q^{-1})$
respectively. Here $s_1^\prime=s_1$, $(s_0^{\prime}f)(x)=f(q^{-1}x^{-1})$
and $\phi_i(x)$ is as in the proof of proposition 5.8, while
$\phi_i^\prime(x)$
is $\phi_i(x)$ with the parameters $(\underline{t},q)$
replaced by $(\underline{t}^{-1},q^{-1})$. In view of the analytic
dependence
on the parameters $\underline{t}$ and $q$, we may assume without loss of
generality
that $0<q<1$ and that the Jordan curve $C$
in the definition 6.2 of $\langle .,. \rangle$ satisfies the following
additional
properties: $C$ has a parametrization of the form $r_C(x)e^{2\pi ix}$
with $r_C: [0,1]\rightarrow (0,\infty)$, and $C$ is $W$-invariant:
$C^{-1}:=\{z^{-1} \, | \, z\in C\}=C$.

For $i=1$ it follows now from the $W$-invariance of $\Delta_+(x)$ that
$\langle T_1f,g\rangle=\langle f, T_1^*g\rangle$ for all $f,g\in
{\mathcal{A}}$
with
\[ T_1^*=k_1-\phi_1(x^{-1})+\frac{\alpha(x)}{\alpha(x^{-1})}\phi_1(x)s_1.\]
Now $\alpha(x)\phi_1(x)=\alpha(x^{-1})\phi_1^\prime(x)$ and
$\phi_1(x^{-1})=\phi_1^\prime(x)$, so that
\[ T_1^*=k_1+\phi_1^\prime(x)(s_1^\prime-1)=(T_1^\prime)^{-1}.
\]
For $i=0$, let $f,g\in {\mathcal{A}}$ and set $h(x)=f(qx^{-1})g(x^{-1})$.
Observe that
\[
(T_0f)(x)g(x^{-1})-f(x)\bigl((T_0^\prime)^{-1}g\bigr)(x^{-1})=\phi_0(x)
\bigl(h(x)-(s_0h)(x)\bigr)
\]
and that
\[ \phi_0(x)\Delta(x)=k_0^{-1}\frac{\bigl(x^2,q^2x^{-2};q\bigr)_{\infty}}
{\bigl(ax,bx,cx,dx,qax^{-1},qbx^{-1},qcx^{-1},qdx^{-1};q\bigr)_{\infty}}
\]
is invariant under $x\mapsto qx^{-1}$.
By the specific properties of $C$, we obtain
\begin{equation*}
\begin{split}
\langle T_0f,g\rangle-\langle f, (T_0^\prime)^{-1}g\rangle&=
\frac{1}{2\pi i}\int_{x\in C}\bigl(h(x)-(s_0h)(x)\bigr)\phi_0(x)
\Delta(x)\frac{dx}{x}\\
&=\frac{1}{2\pi i}\int_{x\in C-qC}h(x)\phi_0(x)\Delta(x)\frac{dx}{x}=0,
\end{split}
\end{equation*}
where the last equality follows from Cauchy's theorem since
$\phi_0(x)\Delta(x)$
is analytic on and within $C-qC$.

The statement for the Dunkl operator $Y$ is immediate
since $Y=T_1T_0$.
\end{proof}

\subsection {}

%6.7
We write $P_m^\prime$ ($m\in {\mathbb{Z}}$)
for the non-symmetric
Askey-Wilson polynomials with respect to the inverse parameters
$(\underline{t}^{-1},q^{-1})$.
\begin{prop}
The two bases $\{P_m\}_{m\in {\mathbb{Z}}}$ and $\{P_n^\prime\}_{n\in
{\mathbb{Z}}}$
of ${\mathcal{A}}$ form a bi-orthogonal system with respect to the
non-degenerate bilinear form $\langle .,. \rangle$, i.e.
$\langle P_m,P_n^\prime\rangle=0$ for $m,n\in {\mathbb{Z}}$ if
$m\not=n$.
\end{prop}
\begin{proof}
By proposition 6.6 and proposition 3.5 we have
\[
\gamma_m\langle P_m,P_n^\prime\rangle=\langle Y\,P_m, P_n^\prime\rangle=
\langle P_m, (Y^\prime)^{-1}P_n'\rangle=\gamma_n\langle
P_m,P_n^\prime\rangle.
\]
It follows that $\langle P_m,P_n^\prime\rangle=0$ if $m\not=n$
since the eigenvalues $\gamma_m$ ($m\in {\mathbb{Z}}$) of $Y$ are
pair-wise different by the conditions 2.11 on the parameters.
\end{proof}

\subsection {}

%6.8
We write $P_m^{+\,\prime}$ and $P_m^{-\,\prime}$ for the symmetric
and anti-symmetric Askey-Wilson polynomial with respect to inverse
parameters $(\underline{t}^{-1},q^{-1})$.

Let ${\mathcal{B}}$ be the basis of ${\mathcal{A}}$ consisting of
$P_m^+$ ($m\in {\mathbb{Z}}_+$) and $P_n^-$ ($n\in {\mathbb{N}}$), and
let ${\mathcal{B}}^\prime$ be the basis of ${\mathcal{A}}$ consisting
of  $P_m^{+\,\prime}$ ($m\in {\mathbb{Z}}_+$)
and $P_m^{-\,\prime}$ ($m\in {\mathbb{N}}$).

\begin{prop}
The pair $({\mathcal{B}}, {\mathcal{B}}^\prime)$
forms a bi-orthogonal system of ${\mathcal{A}}$ with respect to
the non-degenerate bilinear form $\langle .,. \rangle$.
\end{prop}
\begin{proof}
It follows from proposition 6.7 and from the fact that
$P_m^{\pm}\in {\mathcal{A}}(m)=\hbox{span}\{P_m,P_{-m}\}$
for $m\in {\mathbb{Z}}_+$ (with the convention $P_0^-\equiv 0$),
that $\langle P_m^{\pm}, P_n^{\pm\,\prime}\rangle =0$ if
$m\not=n$.

By proposition 6.6 we have
$\langle C_{\pm}f,g\rangle=\langle f, C_{\pm}^\prime g\rangle$ for all
$f,g\in {\mathcal{A}}$, where
$C_{\pm}^\prime=(1+k_1^{\mp 2})^{-1}(1\pm k_1^{\mp 1}T_1^\prime)$ are the
mutually orthogonal, primitive idempotents of $H_0(k_1^{-1})$ (see
5.1) which act on ${\mathcal{A}}$ via $\pi_{\underline{t}^{-1},q^{-1}}$.
Hence $\langle P_m^{\pm}, P_n^{\mp\,\prime}\rangle=0$ for all $m,n\in
{\mathbb{Z}}_+$.
\end{proof}

\subsection {}

%6.9
The bi-orthogonality relations of proposition 6.8 restricted to
${\mathcal{A}}_+={\mathcal{A}}^W$
reduce to the well-known orthogonality relations \cite[theorem 2.3]{AW}
of the symmetric Askey-Wilson polynomials:

\begin{prop}
For all $m\in {\mathbb{Z}}_+$, we have $P_m^{+\,\prime}=P_m^+$.
In particular, $\langle P_m^+,P_n^+\rangle=\bigl(P_m^+,P_n^+\bigr)=0$
if $m\not=n$.
\end{prop}
\begin{proof}
Recall that $P_m^+$ is the unique $W$-invariant Laurent
polynomial of the form $x^m+$ lower order terms with respect to $\preceq$
which is an eigenfunction of the Askey-Wilson $q$-difference operator
$L=(Y+Y^{-1})_{sym}$ with eigenvalue $\gamma_m+\gamma_m^{-1}$.
Then $P_m^+=P_m^{+\,\prime}$ follows from the fact that $L$ and the
eigenvalue $\gamma_m+\gamma_m^{-1}$ are invariant under
$(\underline{t},q)\mapsto (\underline{t}^{-1},q^{-1})$.
The second statement follows now from proposition 6.8 and lemma 6.4.
\end{proof}

%%%%%%%%%%%%%%%%%%%%%%%%%%%%%%%%%%%%%%%%%%%%%%%%%%%%%%%%%%%%%%
%%                                                          %%
%%    The generalized Weyl character formula                %%
%%                                                          %%
%%%%%%%%%%%%%%%%%%%%%%%%%%%%%%%%%%%%%%%%%%%%%%%%%%%%%%%%%%%%%%

\section{The generalized Weyl character formula}

\subsection {}

%7.1
The generalized Weyl character formula relates the
anti-symmetric Askey-Wilson polynomial with the symmetric
Askey-Wilson polynomial via
the generalized Weyl denominator. The generalized Weyl
denominator $\delta(\cdot)$, which we define in the following lemma,
is an explicit anti-symmetric Laurent polynomial of minimal degree
with respect to the total order $\preceq$ on the monomials $\{x^m\}_{m\in
{\mathbb{Z}}}$.

\begin{lem}
We have ${\mathcal{A}}_-=\delta(z)\bigl({\mathcal{A}}_+\bigr)$,
where $\delta=\delta(\,\cdot\,;k_0,k_1)\in {\mathcal{A}}_-$ is given by
\[
\delta(x)=x^{-1}(x-k_0^{-1}k_1^{-1})(x+k_0k_1^{-1})=x^{-1}(x-a^{-1})
(x-b^{-1}).
\]
\end{lem}
\begin{proof}
By 2.21 we have $(T_1+k_1^{-1})\delta(z)=\delta(z^{-1})(T_1-k_1)$.
Combined with 5.2 this implies
$\delta(z)\bigl({\mathcal{A}}_+\bigr) \subseteq {\mathcal{A}}_-$.
Let now $f\in {\mathcal{A}}_-$, and set $g=\delta^{-1}f\in
{\mathbb{C}}(x)$.
Using the extended action of $T_1$ on ${\mathbb{C}}(x)$, see 5.6,
we derive that
\[
\delta(x^{-1})\bigl((T_1-k_1)g\bigr)(x)=\bigl((T_1+k_1^{-1})f\bigr)(x)=0,
\]
so that $\bigl((T_1-k_1)g\bigr)(x)=0$. Since $T_1-k_1=\phi_1(x)(s_1-1)$
with
$0\not=\phi_1(x)\in {\mathbb{C}}(x)$, we conclude that $g$ is $W$-invariant
in ${\mathbb{C}}(x)$.
In particular,  $\delta(x^{-1})f(x)=\delta(x)f(x^{-1})$ in ${\mathcal{A}}$.
Since $\delta(x^{-1})=a^{-1}b^{-1}x^{-1}(x-a)(x-b)$ and $\delta(x)$
are relative coprime in the unique factorisation domain
${\mathcal{A}}$ by the conditions 2.11 on the parameters,
we conclude that $f$ is divisible by $\delta$ in ${\mathcal{A}}$.
Hence, $g=\delta^{-1}f\in {\mathcal{A}}$. Since $g$ is furthermore
$W$-invariant, we
conclude that $g\in {\mathcal{A}}_+$, hence
${\mathcal{A}}_-\subseteq \delta(z)\bigl({\mathcal{A}}_+\bigr)$.
\end{proof}

\subsection {}

%7.2
The bilinear form $\langle .,. \rangle$ restricted to
${\mathcal{A}}_-$
can now be related to the bilinear
form $\bigl( .,. \bigr)$ on ${\mathcal{A}}_+$ using lemma 7.1.
We identify $\underline{t}$ with $(k_0,k_1,u_0,u_1)$ in accordance with
2.11.

\begin{lem}
Assume that the parameters satisfy the additional conditions 6.1.
Let $\delta^\prime(x)=\delta(x;k_0^{-1},k_1^{-1})$. Then for all $f,g\in
{\mathcal{A}}^W$,
\[\langle \delta(z)f, \delta^\prime(z)g\rangle_{\underline{t},q}=
\frac{1}{2}(1+k_1^{-2})\bigl(f,g\bigr)_{k_0,qk_1,u_0,u_1,q}.
\]
\end{lem}
\begin{proof}
Set $\alpha^\prime(x)=\alpha(x;k_1^{-1},u_1^{-1})$, see 6.2, then
\[\delta(x)\delta^\prime(x^{-1})\alpha(x)=
(1-ax)(1-bx)(1-ax^{-1})(1-bx^{-1})\alpha^\prime(x).
\]
By the explicit expression for the $W$-invariant weight function
$\Delta_+(x;\underline{t};q)$, see 6.2, we obtain
\[\langle \delta(z)f,\delta^\prime(z)g\rangle_{\underline{t},q}=
\frac{1}{2\pi i}\int_{x\in C}
f(x)g(x^{-1})\alpha^\prime(x)\Delta_+(x;k_0,qk_1,u_0,u_1;q)\frac{dx}{x}
\]
for $f,g\in {\mathcal{A}}^W$.
The result follows now by symmetrizing the integrand, cf. 6.4.
\end{proof}

\subsection {}

%7.3
We are now in a position to relate the
anti-symmetric Askey-Wilson polynomial $P_m^-$ ($m\in
{\mathbb{N}}$) with the symmetric Askey-Wilson polynomial
$P_{m-1}^+$ via the generalized Weyl denominator $\delta$.
The result is as follows.

\begin{prop}[Generalized Weyl character formula]
For $m\in {\mathbb{N}}$ we have
\[
P_m^-(x;\underline{t};q)=\delta(x;k_0,k_1)P_{m-1}^+(x;k_0,qk_1,u_0,u_1;q).
\]
\end{prop}
\begin{proof}
We first prove the proposition when $|q|<1$.

We assume for the
moment that the multiplicity function $\underline{t}$
satisfies the additional conditions 6.1 and that
$\bigl( .,. \bigr)_{k_0,qk_1,u_0,u_1,q}$ restricts to a
non-degenerate bilinear form on ${\mathcal{A}}_m^W:=
\hbox{span}\{m_n \, | \, n=0,\ldots,m-1\}$, where
$m_0(x)=1$ and $m_n(x)=x^n+x^{-n}$ for $n\in {\mathbb{N}}$.
These are generic conditions on the parameters, which can be removed by
continuity
at the end of the proof. Indeed, observe that the restriction of
$\bigl( .,. \bigr)_{\underline{t},q}$ to ${\mathcal{A}}_m^W$
is non-degenerate when $0<a,b,c,d,q<1$,
since then the bilinear form can be given as integration over the unit
circle with respect to the positive weight function $\Delta_+(x)$.
By analytic continuation, it follows that the restriction of
$\bigl( .,. \bigr)_{\underline{t},q}$ to ${\mathcal{A}}_m^W$
is non-degenerate for generic parameter values
satisfying the conditions 6.1.

By proposition 6.9 we conclude that
$P_{m-1}^+(x;k_0,qk_1,u_0,u_1;q)$ is the unique $W$-invariant
Laurent polynomial of the form $x^{m-1}+$ lower order terms w.r.t.
$\preceq$ which is orthogonal to $m_n$ for $n=0,\ldots,m-2$ with
respect to the bilinear form $\bigl( .,.
\bigr)_{k_0,qk_1,u_0,u_1,q}$. We show that
$p(x)=\delta(x;k_0,k_1)^{-1}P_m^-(x;\underline{t};q)$
satisfies the same characterizing
conditions.

By lemma 7.1 we have $p\in {\mathcal{A}}^W$.
Since $P_m^-(x)=x^m+$ lower order terms w.r.t. $\preceq$
and $\delta(x)=x+$ lower order terms w.r.t $\preceq$, we have
$p(x)=x^{m-1}+$ lower order terms w.r.t. $\preceq$.
By the triangularity properties of the anti-symmetric Askey-Wilson
polynomials (see proposition 5.3{\bf (ii)}) and by lemma 7.1, we see that
$\delta^\prime(z)m_n\in \hbox{span}\{P_k^{-\,\prime}\, | \,
k=1,\ldots,m-1\}$
for $n=0,\ldots,m-2$. By lemma 7.2 and proposition 6.8 we conclude that
\[
\frac{1}{2}(1+k_1^{-2})\bigl(p,m_n\bigr)_{k_0,qk_1,u_0,u_1,q}=
\langle P_m^-,\delta^\prime(z)m_n\rangle_{\underline{t},q}=0,\quad
(n=0,\ldots,m-2).
\]
Hence $p(x)=P_{m-1}^+(x;k_0,qk_1,u_0,u_1;q)$, as desired.

The proof for $|q|>1$ is a slight
modification of the arguments for $|q|<1$. One uses now that
$P_m^{+\,\prime}=P_m^+$ for $m\in {\mathbb{Z}}_+$ (see proposition
6.9) and a characterization of $P_m^+(x;\underline{t};q)$
in terms of the bilinear form $\bigl( .,.\bigr)_{\underline{t}^{-1},q^{-1}}$.
\end{proof}

\subsection {}

%7.4
We have now all the ingredients to express the diagonal terms
$\langle P_m, P_m^\prime\rangle$ ($m\in {\mathbb{Z}}$) and
$\langle P_m^-, P_m^{-\,\prime}\rangle$ ($m\in {\mathbb{N}}$)
of the bi-orthogonality relations in proposition 6.7 and
proposition 6.8 in terms of the ``quadratic norms''
$\langle P_m^+, P_m^+\rangle=\frac{1}{2}(1+k_1^2)\bigl(P_m^+,
P_m^+\bigr)$ ($m\in {\mathbb{Z}}_+$).

Indeed, by the generalized Weyl character formula, lemma
7.2 and proposition 6.9 we have for $m\in {\mathbb{N}}$
that
\begin{equation*}
\begin{split}
&\langle P_m^-,
P_m^{-\,\prime}\rangle_{\underline{t},q}=\\
&\,\,\frac{1}{2}(1+k_1^{-2})\bigl(P_{m-1}^+(\,\cdot\,;k_0,qk_1,u_0,u_1;q),
P_{m-1}^+(\,\cdot\,;k_0,qk_1,u_0,u_1;q)\bigr)_{k_0,qk_1,u_0,u_1,q}.
\end{split}
\end{equation*}
On the other hand, for $m\in {\mathbb{N}}$ we have
$P_m^-=-(1+k_1^{-2})C_-P_{-m}$ (cf. the proof of proposition 5.3),
so that
\[ \langle P_{-m},
P_{-m}^\prime\rangle=\frac{(1+k_1^2)(1-\gamma_m^{-2})}
{(1+k_0k_1^{-1}\gamma_m^{-1})(1-k_0^{-1}k_1^{-1}\gamma_m^{-1})}
\langle P_m^-, P_m^{-\,\prime}\rangle
\]
by proposition 5.3{\bf (ii)}. Here we have used that $\langle
C_+f,g\rangle=
\langle f,C_-^\prime g\rangle$ for all $f,g\in {\mathcal{A}}$
(see the proof of proposition 6.8) and that $C_-P_m^-=P_m^-$.
Similarly, we can relate $\langle P_m, P_m^\prime\rangle$ for $m\in
{\mathbb{N}}$
to $\langle P_{m}^-,P_{m}^{-\,\prime}\rangle$ using (the proof of)
proposition 5.3{\bf (ii)}.

\subsection {}

%7.5
The generalized Weyl character formula plays a crucial role in
the study of shift operators for the symmetric Askey-Wilson
polynomials. In turn, shift operators can be used to explicitly evaluate
the quadratic norms $\bigl(P_m^+, P_m^+\bigr)_{\underline{t},q}$ ($m\in
{\mathbb{Z}}_+$). Combined with 7.4, this leads to explicit
evaluations of all the diagonal terms of
the bi-orthogonality relations in proposition 6.7 and proposition 6.8.

In section 11 we present another method for deriving explicit
expressions of the diagonal terms, which uses the double affine Hecke
algebra in an essential way. This method gives
more insight in the particular structure of the diagonal
terms. Namely, it shows that the diagonal terms
can be naturally expressed in terms of the
residue of the weight function in a certain
specific simple
pole, the constant term $\langle 1,1\rangle$, and the value of
the Askey-Wilson polynomial at the point
$a^{-1}=k_1^{-1}u_1^{-1}$.

We return to shift operators in section 12 in order to evaluate
the constant term $\langle 1,1\rangle$ (which is the well-known
Askey-Wilson integral, see \cite{AW}).

%%%%%%%%%%%%%%%%%%%%%%%%%%%%%%%%%%%%%%%%%%%%%%%%%%%%%%%%%%%%%%%%%%%%%
%%                                                                 %%
%%               The double affine Hecke algebra                   %%
%%                                                                 %%
%%%%%%%%%%%%%%%%%%%%%%%%%%%%%%%%%%%%%%%%%%%%%%%%%%%%%%%%%%%%%%%%%%%%%

\section{The double affine Hecke algebra}
\subsection {}

%8.1
Recall from 2.21 that the double affine Hecke algebra
${\mathcal{H}}(S;\underline{t};q)$ is the
subalgebra of $\hbox{End}_{\mathbb{C}}\bigl({\mathcal{A}}\bigr)$
generated by $\pi_{\underline{t},q}\bigl(H(R;k_0,k_1)\bigr)$
and by ${\mathcal{A}}$, where ${\mathcal{A}}$ is regarded as
subalgebra of $\hbox{End}_{\mathbb{C}}({\mathcal{A}})$ via its regular
representation.

\subsection {}

%8.2
We have observed in 2.22
that there is a unique surjective algebra homomorphism
$\phi: {\mathcal{F}}(\underline{t};q)\rightarrow
{\mathcal{H}}(S;\underline{t};q)$ satisfying the conditions
as stated in theorem 2.22. In particular, we have
\[z^{a_0^\vee}T_0=T_0^{-1}z^{-a_0^\vee}+u_0^{-1}-u_0,\qquad
z^{-a_1^\vee}T_1^{-1}=T_1z^{a_1^\vee}+u_1-u_1^{-1}
\]
in ${\mathcal{H}}(S;\underline{t};q)$.
Combined with 2.15, we conclude that $\mathcal{H}(S;\underline{t};q)$
is spanned by the elements $z^mY^n$ and $z^mT_1Y^n$
where $m,n\in {\mathbb{Z}}$.
In fact, we have the following stronger result, see Sahi
\cite[theorem 3.2]{S} for the higher rank setting.

\begin{prop}
The set $\{z^mY^n, z^mT_1Y^n \, | \, m,n\in {\mathbb{Z}} \}$ is
a linear basis of $\mathcal{H}(S;\underline{t};q)$.
\end{prop}
\begin{proof}
Assume that
$X=\sum_{m,n\in
{\mathbb{Z}}}\bigl(c_{m,n}^1z^mY^n+c_{m,n}^2z^mT_1Y^n\bigr)=0$
in $\hbox{End}_{\mathbb{C}}\bigl({\mathcal{A}}\bigr)$
with only finitely many
coefficients
$c_{m,n}^j$ non-zero. Since $z$ is invertible
in $\hbox{End}_{\mathbb{C}}\bigl({\mathcal{A}}\bigr)$,
we may assume without loss of
generality that $c_{m,n}^j=0$ unless $m\in {\mathbb{Z}}_+$. Suppose that
not all coefficients $c_{m,n}^j$ are zero.
Let $m_0\in {\mathbb{Z}}_+$ be the
largest positive integer such that $c_{m_0,n}^j$ is non-zero for some
$n\in {\mathbb{Z}}$ and some $j\in \{1,2\}$.

Let $X$ act on the non-symmetric Askey-Wilson polynomial $P_{-l}$
($l\in {\mathbb{N}}$),
and consider the coefficient of $x^{m_0+l}$
in the resulting expression using proposition 3.5 and proposition 4.1.
We obtain $\sum_{n\in {\mathbb{Z}}_+}k_1c_{m_0,n}^2\gamma_l^{-n}=0$ for all
$l\in {\mathbb{N}}$, hence $c_{m_0,n}^2=0$ for all $n\in {\mathbb{Z}}$.

Let now $X$ act on $P_l$ with $l\in {\mathbb{Z}}_+$, and again consider the
coefficient of $x^{m_0+l}$ in the resulting expression. Then
$\sum_{n\in {\mathbb{Z}}_+}c_{m_0,n}^1\gamma_l^n=0$
for all $l\in {\mathbb{Z}}_+$,
hence $c_{m_0,n}^1=0$ for all
$n\in {\mathbb{Z}}$. This gives the desired contradiction.
\end{proof}

\subsection {}

%8.3
We can finish now the proof of theorem 2.22 using proposition 8.2.
It suffices to show that the surjective algebra homomorphism
$\phi: {\mathcal{F}}(\underline{t};q)\rightarrow
\mathcal{H}(S;\underline{t};q)$
defined in theorem 2.22 is injective. Set
$w:=V_1^{-1}(V_1^\vee)^{-1}=
q^{1/2}V_0V_0^\vee\in {\mathcal{F}}(\underline{t};q)$,
then $\phi(w)=z$. Observe that
$V_0^\vee=q^{1/2}w^{-1}V_0+u_0-u_0^{-1}$
and that
$V_1^\vee=w^{-1}V_1^{-1}$ in ${\mathcal{F}}(\underline{t};q)$, so that
${\mathcal{F}}(\underline{t};q)$ is generated by $w^{\pm 1}, V_0$ and
$V_1$ as an algebra.
Furthermore,
\[ V_0w=qw^{-1}V_0+(k_0-k_0^{-1})w+q^{1/2}(u_0-u_0^{-1}),
\qquad
V_1w=w^{-1}V_1^{-1}+u_1^{-1}-u_1
\]
in ${\mathcal{F}}(\underline{t};q)$, so that
any element in ${\mathcal{F}}(\underline{t};q)$
can be written as a finite linear combination of elements of the form
$f(w)X$, where $f(w)$ is a Laurent polynomial in $w$ and $X$ is an element
in the subalgebra of ${\mathcal{F}}(\underline{t};q)$ generated by
$V_0$ and $V_1$. By 2.15 and by
the relations (1) in theorem 2.22 for the generators
$V_0,V_1\in {\mathcal{F}}(\underline{t};q)$
it follows that ${\mathcal{F}}(\underline{t};q)$ is spanned
by
$\{w^mZ^n, w^mV_1Z^n \, | \, m,n\in {\mathbb{Z}} \}$, where $Z:=V_1V_0$.
Since the image of these elements under $\phi$ are linear independent
by proposition 8.2, we conclude that $\phi$ is injective.

\subsection {}

%8.4
In the remainder of the paper we use the notations
$T_0^\vee:= T_0^{-1}z^{-a_0^\vee}\in
\mathcal{H}(S;\underline{t},q)$
and $T_1^\vee:= z^{-a_1^\vee}T_1^{-1}\in \mathcal{H}(S;\underline{t};q)$
for the images of $V_0^\vee$ and $V_1^\vee$ respectively under the algebra
isomorphism $\phi: {\mathcal{F}}(\underline{t};q)\rightarrow
{\mathcal{H}}(S;\underline{t};q)$ (see theorem 2.22).

\subsection {}

%8.5
We associate with the multiplicity function $\underline{t}\simeq
(k_0,k_1,u_0,u_1)$
a dual multiplicity function $\underline{\tilde{t}}$ by interchanging $k_0$
and $u_1$, so $\underline{\tilde{t}}\simeq (u_1,k_1,u_0,k_0)$.
We write $\widetilde{T}_0,\widetilde{T}_1, \widetilde{T}_0^\vee,
\widetilde{T}_1^\vee$
for the generators of $\mathcal{H}(S;\underline{\tilde{t}};q)$
(cf. 2.22 and 8.4),
and we write
$\widetilde{Y}=\widetilde{T}_1\widetilde{T}_0$ for the associated Dunkl
operator
and $\widetilde{z}=\widetilde{T}_1^{-1}(\widetilde{T}_1^\vee)^{-1}=
q^{1/2}\widetilde{T}_0\widetilde{T}_0^\vee$ for the corresponding
``multiplication by $x$'' operator. The first part of the
following proposition is a special case of Sahi's results
in \cite[section 7]{S}.

\begin{prop}
{\bf (i)}
The application $T_0\mapsto \widetilde{T}_1^\vee$,
$T_1\mapsto \widetilde{T}_1$, $T_0^\vee\mapsto \widetilde{T}_0^\vee$ and
$T_1^\vee\mapsto
\widetilde{T}_0$ uniquely extend to an anti-algebra isomorphism
$\nu=\nu_{\underline{t},q}: \mathcal{H}(S;\underline{t};q)\rightarrow
\mathcal{H}(S;\underline{\tilde{t}};q)$. Furthermore,
$\nu_{\underline{t},q}^{-1}=\nu_{\underline{\tilde{t}},q}$ and
$\nu(z)=\widetilde{Y}^{-1}$, $\nu(Y)=\widetilde{z}^{-1}$.

{\bf (ii)} The application $T_0\mapsto
\widetilde{T}_1^{-1}\widetilde{T}_1^\vee\widetilde{T}_1$,
$T_1\mapsto \widetilde{T}_1$, $T_0^\vee\mapsto
\widetilde{T}_0\widetilde{T}_0^\vee
\widetilde{T}_0^{-1}$ and $T_1^\vee\mapsto \widetilde{T}_0$ uniquely extend
to an
algebra isomorphism
$\mu=\mu_{\underline{t},q}: \mathcal{H}(S;\underline{t};q)\rightarrow
\mathcal{H}(S;\underline{\tilde{t}};q)$. Furthermore,
$\mu(Y)=\widetilde{z}^{-1}$.
\end{prop}

\begin{proof}
By theorem 2.22 it suffices to check that $\mu$ (respectively $\nu$)
is compatible with the defining relations in
${\mathcal{F}}(\underline{t};q)\simeq \mathcal{H}(S;\underline{t};q)$.
This can be done by direct computations. It is immediate that
$\nu_{\underline{\tilde{t}},q}$ is the inverse of $\nu_{\underline{t},q}$.

Observe that the application $\widetilde{T}_0\rightarrow T_1^\vee$,
$\widetilde{T}_1\mapsto T_1$, $\widetilde{T}_0^\vee\mapsto
(T_1^\vee)^{-1}T_0^\vee T_1^\vee$ and $\widetilde{T}_1^\vee\mapsto
T_1T_0T_1^{-1}$
uniquely extend to an algebra homomorphism from
$\mathcal{H}(S;\underline{\tilde{t}};q)$ to
$\mathcal{H}(S;\underline{t};q)$. It is immediate that this
homomorphism is the inverse of $\mu$.
\end{proof}

\subsection {}

%8.6
Following the terminology of Sahi \cite[section 7]{S}, we call
$\nu=\nu_{\underline{t},q}$
the  {\it duality anti-isomorphism}. Furthermore, we call
$\mu=\mu_{\underline{t},q}$
the {\it duality isomorphism}.  These duality isomorphisms
play a fundamental role in the theory of non-symmetric Askey-Wilson
polynomials.
In particular, the duality anti-isomorphism can be used to show that the
geometric parameter $x$ and the spectral parameter $\gamma$ of the
non-symmetric Askey-Wilson polynomial are in a sense interchangeable
(see Sahi \cite[section 7]{S} or section 10).
The duality isomorphism describes the
intertwining properties of the action of the double affine Hecke
algebra under the non-symmetric Askey-Wilson transform, see section 11.

\subsection {}

%8.7
We write $T_i^\prime$ and $T_i^{\vee \prime}$ \textup{(}$i=0,1$\textup{)}
for the generators of $\mathcal{H}(S;\underline{t}^{-1};q^{-1})$,
cf. 2.22, 6.6 and 8.4.

\begin{prop}
There exists a unique anti-algebra isomorphism ${}^*:
{\mathcal{H}}(S;\underline{t};q)
\rightarrow {\mathcal{H}}(S;\underline{t}^{-1};q^{-1})$ such that
$T_i^* =(T_i^\prime)^{-1}$ and
$(T_i^\vee)^*=(T_i^{\vee \prime})^{-1}$ for $i=0,1$. Furthermore, $T^*$
coincides with the adjoint of $T\in {\mathcal{H}}(S;\underline{t};q)$
if the parameters satisfy the additional conditions 6.1.
\end{prop}
\begin{proof}
The first statement follows easily from theorem 2.22.
For the second statement, it suffices to compute the adjoint of $T_i^\vee$
($i=0,1$) in view of proposition 6.6.
Let $z^\prime=q^{-1/2}T_0^\prime T_0^{\vee \prime}$ be the ``multiplication
by $x$''
operator in $\mathcal{H}(S;\underline{t}^{-1};q^{-1})$.
It is immediate that  $z^*=(z^\prime)^{-1}$. Combined with proposition
6.6 we obtain $(T_0^\vee)^*=q^{-1/2}(z^\prime)^{-1}T_0^\prime=
(T_0^{\vee \prime})^{-1}$
and $(T_1^\vee)^*=T_1^\prime z^\prime=(T_1^{\vee \prime})^{-1}$.
This gives the desired result.
\end{proof}

%%%%%%%%%%%%%%%%%%%%%%%%%%%%%%%%%%%%%%%%%%%%%%%%%%%%%%%%%%%%%%%%%%
%%                                                              %%
%%      Intertwiners as creation operators                      %%
%%                                                              %%
%%%%%%%%%%%%%%%%%%%%%%%%%%%%%%%%%%%%%%%%%%%%%%%%%%%%%%%%%%%%%%%%%%

\section{Intertwiners as creation operators}

\subsection {}

In corollary 2.17 we have introduced the non-affine intertwiner $S_1$
and derived its basic property.
The results of the previous section allow us to derive the
following analogous result for the commutator
$[Y, T_1^\vee]\in {\mathcal{H}}(S;\underline{t};q)$, see Sahi
\cite[theorem 5.1]{S} for the result in the higher rank setting.
\begin{cor}[The affine intertwiner]
Set $S_0:=[Y, T_1^\vee]=YT_1^\vee-T_1^\vee Y\in
{\mathcal{H}}(S;\underline{t};q)$.
Then $g(Y)S_0=S_0g(q^{-1}Y^{-1})$
for all $g(Y)\in {\mathbb{C}}[Y^{\pm 1}]$.
\end{cor}
\begin{proof}
By 2.21 we have $f(z)[T_0,z^{-1}]=[T_0,z^{-1}](s_0f)(z)$ in
$\mathcal{H}(S;\underline{t};q)$ for any Laurent polynomial $f$.
Apply now the duality anti-isomorphism $\nu_{\underline{t},q}$
to this equality, and replace the parameters by dual parameters in the
resulting identity. This gives $S_0f(Y^{-1})=
(s_0f)(Y^{-1})S_0$ for any Laurent polynomial $f$. The corollary is
now immediate.
\end{proof}

\subsection {}

%9.2
Recall from corollary 4.3 that the action of $S_1$ on the non-symmetric
Askey-Wilson polynomials is completely explicit. In the following
lemma we give the analogous result
for the action of the affine intertwiner $S_0$ on $P_m$ ($m\in
{\mathbb{Z}}_+$).

\begin{lem}
Let $m\in {\mathbb{Z}}_+$, then
$S_0(P_m)=\bigl(\gamma_{-m-1}-\gamma_m\bigr)k_1^{-1}P_{-m-1}$.
\end{lem}
\begin{proof}
Let $m\in {\mathbb{Z}}_+$. By proposition 4.1 we have
\[S_0(P_m)=(Y-\gamma_m)\bigl((\alpha_m+k_1^{-1}-k_1)
z^{-1}P_m+\beta_mz^{-1}P_{-m}\bigr).
\]
It follows then from proposition 3.4 that the leading term of
$S_0(P_m)$ with respect to the total order $\preceq$ on the
monomials equals
\[(\gamma_{-m-1}-\gamma_{m})\bigl((\alpha_m+
k_1^{-1}-k_1)c_m+\beta_m\bigr)x^{-m-1},
\]
where $c_0:=1$ and where $c_m$ ($m\in {\mathbb{N}}$) is the unique
constant such
that $P_m(x)=x^m+c_mx^{-m}+$ lower order terms
w.r.t $\preceq$.
By proposition 5.3{\bf (i)} we have
\[c_m=1-\frac{(1+k_0k_1^{-1}
\gamma_m)(1-k_0^{-1}k_1^{-1}\gamma_m)}{(1-\gamma_m^2)}
=k_1^{-1}\alpha_m.
\]
Furthermore, recall from the proof of proposition 4.1 that
$\beta_m=k_1^{-1}(k_1-\alpha_m)(k_1^{-1}+\alpha_m)$. Hence
the leading term of $S_0(P_m)$ reduces
to $(\gamma_{m-1}-\gamma_m)k_1^{-1}x^{-m-1}$.
On the other hand, corollary 9.1 implies that $S_0(P_m)=d_mP_{-m-1}$
for some constant $d_m$. By the leading term considerations, we
conclude that $d_m=(\gamma_{m-1}-\gamma_m)k_1^{-1}$.
\end{proof}

\subsection {}

%9.3

The intertwiners $S_0$ and $S_1$
can be used to create the non-symmetric Askey-Wilson
polynomial $P_m$ ($m\in {\mathbb{Z}}$)
from the unit polynomial $1\in {\mathcal{A}}$
in the following way.

\begin{prop}
We have $(S_1S_0)^m(1)=d_mP_m$ for $m\in {\mathbb{Z}}_+$
and $\bigl(S_0(S_1S_0)^{m-1}\bigr)(1)=d_{-m}P_{-m}$ for $m\in
{\mathbb{N}}$,
with the constants $d_m$ \textup{(}$m\in {\mathbb{Z}}$\textup{)} given by
\begin{equation*}
\begin{split}
d_m&=q^{-(m+1)m}k_0^{-2m}k_1^{-2m}\bigl(qk_0^2k_1^2;q\bigr)_{2m},
\qquad m\in {\mathbb{Z}}_+,\\
d_{-m}&=q^{-m^2}k_0^{1-2m}k_1^{-2m}\bigl(qk_0^2k_1^2;q\bigr)_{2m-1},
\qquad m\in {\mathbb{N}}.
\end{split}
\end{equation*}
\end{prop}

\begin{proof}
By corollary 2.17 and corollary 9.1
we have $g(Y)(S_1S_0)=(S_1S_0)g(qY)$ for all
$g\in {\mathcal{A}}$.
It follows that $F_m:=(S_1S_0)^m(1)\in {\mathcal{A}}$ for $m\in
{\mathbb{Z}}_+$
satisfies $g(Y)F_m=g(\gamma_m)F_m$ for all $m\in {\mathbb{Z}}$, so
$F_m=d_mP_m$ for some constant $d_m$.
Similarly, we obtain $F_{-m}:=S_0(S_1S_0)^{m-1}(1)=d_{-m}P_{-m}$ for some
constant $d_{-m}$ when $m\in {\mathbb{N}}$.
By corollary 4.3 and lemma 9.2, we have the recurrence relations
\begin{equation*}
\begin{split}
d_m&=(\gamma_{-m}-\gamma_m)k_1d_{-m},\qquad m\in {\mathbb{N}},\\
d_{-m-1}&=(\gamma_{-m-1}-\gamma_m)k_1^{-1}d_m,\qquad m\in
{\mathbb{Z}}_+.
\end{split}
\end{equation*}
Together with the initial condition $d_0=1$, we obtain the explicit
expressions for $d_m$ ($m\in {\mathbb{Z}}$)
by complete induction with respect to $m$.
\end{proof}

%%%%%%%%%%%%%%%%%%%%%%%%%%%%%%%%%%%%%%%%%%%%%%%%%%%%%%%%%%%%%%%%%%%%
%%                                                                %%
%%     Evaluation formula and duality                             %%
%%                                                                %%
%%%%%%%%%%%%%%%%%%%%%%%%%%%%%%%%%%%%%%%%%%%%%%%%%%%%%%%%%%%%%%%%%%%%

\section{Evaluation formula and duality}

\subsection {}

%10.1
Let $\hbox{Ev}=\hbox{Ev}_{\underline{t},q}:
\mathcal{H}(S;\underline{t};q)\rightarrow {\mathbb{C}}$
be the linear map defined by
$\hbox{Ev}(X):=\bigl(X(1)\bigr)(k_1^{-1}u_1^{-1})$, where
$1\in {\mathcal{A}}$ is the Laurent polynomial identically equal to one.
Observe that $\hbox{Ev}$ satisfies
\[\hbox{Ev}\bigl(T_1^{\pm 1}X\bigr)=k_1^{\pm 1}\hbox{Ev}(X),\qquad
X\in \mathcal{H}(S;\underline{t};q),
\]
since $(T_1f)(k_1^{-1}u_1^{-1})=k_1f(k_1^{-1}u_1^{-1})$ for all $f\in
{\mathcal{A}}$
by the explicit expression 2.19
for the difference-reflection operator $T_1$.

\subsection {}

%10.2
Observe that we can evaluate $\hbox{Ev}(P_m(z))=P_m(k_1^{-1}u_1^{-1})$
explicitly using proposition 5.10, since the two ${}_4\phi_3$'s
in the right hand side of the formula for $P_m(x)$ are
equal to one when $x=a^{-1}=k_1^{-1}u_1^{-1}$.
We give here an alternative, inductive proof for the evaluation which
only uses the Rodrigues type formula for the non-symmetric
Askey-Wilson polynomials
in terms of the intertwiners $S_0$ and $S_1$, see proposition 9.3.
We abuse notation by writing
$\hbox{Ev}(f)=\hbox{Ev}(f(z))=f(k_1^{-1}u_1^{-1})$ for $f\in
{\mathcal{A}}$.

\begin{prop}
{\bf (i)} For $m\in {\mathbb{Z}}_+$, we have
\begin{equation*}
\begin{split}
\hbox{Ev}\bigl(P_m\bigr)&=k_1^{-m}u_1^{-m}
\frac{\bigl(-qk_1^2,q^{1/2}k_0k_1u_0u_1,
-q^{1/2}k_0k_1u_0^{-1}u_1;q\bigr)_m}
{\bigl(q^{m+1}k_0^2k_1^2;q\bigr)_m}\\
&=a^{-m}\frac{\bigl(qab,ac,ad;q\bigr)_m}
{\bigl(q^mabcd;q\bigr)_m}.
\end{split}
\end{equation*}
{\bf (ii)} For $m\in {\mathbb{N}}$, we have
\begin{equation*}
\begin{split}
\hbox{Ev}\bigl(P_{-m}\bigr)&=
\frac{k_1^{-m}u_1^{-m}}{1+k_1^{-2}}\frac{\bigl(-k_1^2,q^{1/2}k_0k_1u_0u_1,
-q^{1/2}k_0k_1u_0^{-1}u_1;q\bigr)_m}{\bigl(q^{m}k_0^2k_1^2;q\bigr)_m}\\
&=
\frac{a^{-m}}{(1-a^{-1}b^{-1})}\frac{\bigl(ab,ac,ad;q\bigr)_m}
{\bigl(q^{m-1}abcd;q\bigr)_m}.
\end{split}
\end{equation*}
\end{prop}
\begin{proof}
We write $F_m=(S_1S_0)^m(1)$ for $m\in {\mathbb{Z}}_+$ and
$F_{-m}=(S_0(S_1S_0)^{m-1})(1)$ for $m\in {\mathbb{N}}$,
so that $F_m=d_mP_m$ for $m\in {\mathbb{Z}}$
with the specific constants $d_m$ as given in proposition 9.3.
The proposition follows then from the explicit evaluation of the
$d_m$, see proposition 9.3, and from the recurrence
relations
\[
\hbox{Ev}(F_m)=k_1^{-1}\gamma_{-m}(1-k_0k_1\gamma_m)(1+k_0^{-1}k_1\gamma_m)
\hbox{Ev}(F_{-m}),
\qquad m\in {\mathbb{N}},
\]
respectively
\[
\hbox{Ev}(F_{-m})=u_1^{-1}\gamma_{-m}(1-u_0u_1q^{-1/2}\gamma_m)
(1+u_0^{-1}u_1q^{-1/2}\gamma_m)\hbox{Ev}(F_{m-1}),\quad
m\in {\mathbb{N}},
\]
by complete induction with respect to $m$.
Let $m\in {\mathbb{N}}$. For the first recurrence relation, observe
that by formula 10.1
and by $YF_{-m}=\gamma_{-m}F_{-m}$ we have
\[\hbox{Ev}(F_m)=\hbox{Ev}(S_1F_{-m})=
\hbox{Ev}\bigl((k_1\gamma_{-m}-k_1T_0T_1)F_{-m}\bigr).
\]
To reduce the $T_0T_1$-term, we use the relation
\[
T_0T_1=Y^{-1}+(k_0-k_0^{-1})T_1+(k_1-k_1^{-1})T_1^{-1}Y-
(k_0-k_0^{-1})(k_1-k_1^{-1})
\]
in ${\mathcal{H}}$ and formula 10.1, which yields
\[\hbox{Ev}(F_m)=k_1^{-1}\gamma_{-m}(1-k_0k_1\gamma_m)
(1+k_0^{-1}k_1\gamma_m)
\hbox{Ev}(F_{-m})
\]
after a direct computation.
For the second recurrence relation,
observe that
\[\hbox{Ev}(F_{-m})=\hbox{Ev}(S_0F_{m-1})=
\hbox{Ev}\bigl((Yz^{-1}T_1^{-1}-u_1\gamma_{m-1})F_{m-1}\bigr)
\]
by formula 10.1, since $YF_{m-1}=\gamma_{m-1}F_{m-1}$.
To reduce the $Yz^{-1}T_1^{-1}$-term, we use the relation
\[
Yz^{-1}T_1^{-1}=q^{-1}z^{-1}T_1^{-1}Y^{-1}+
q^{-1}(u_1^{-1}-u_1)Y^{-1}+q^{-1/2}(u_0^{-1}-u_0)
\]
in ${\mathcal{H}}$ and formula 10.1, which gives the desired recursion
\[
\hbox{Ev}(F_{-m})=u_1^{-1}\gamma_{-m}(1-u_0u_1q^{-1/2}\gamma_m)
(1+u_0^{-1}u_1q^{-1/2}\gamma_m)\hbox{Ev}(F_{m-1})
\]
after a direct computation.
This completes the proof of the proposition.
\end{proof}

\subsection {}

%10.3
The explicit expression for the
(anti-)symmetric Askey-Wilson polynomial as linear combination of
non-symmetric Askey-Wilson polynomials (see proposition 5.3) can be used to
express $\hbox{Ev}(P_m^{\pm})$ as a linear combination of
$\hbox{Ev}(P_m)$ and $\hbox{Ev}(P_{-m})$. Combined with
proposition 10.2, this leads to explicit evaluation formulas
for the (anti-)symmetric Askey-Wilson polynomials $P_m^{\pm}$. In
particular, it follows that
\[\hbox{Ev}(P_m^+)=
P_m^+(a)=a^{-m}\frac{\bigl(ab,ac,ad;q\bigr)_m}
{\bigl(q^{m-1}abcd;q\bigr)_m},\qquad m\in {\mathbb{Z}}_+.
\]
This result can also be obtained directly from the explicit
expression of $P_m^+$ in terms of a terminating, balanced
${}_4\phi_3$, see theorem 5.9.

\subsection {}

%10.4
The evaluation mapping $\hbox{Ev}$ and the duality anti-isomorphism $\nu$
are compatible in the following sense.

\begin{lem}
For all $X\in \mathcal{H}(S;\underline{t};q)$ we have
$\hbox{Ev}_{\underline{\tilde{t}},q}\bigl(\nu_{\underline{t},q}(X)\bigr)=
\hbox{Ev}_{\underline{t},q}\bigl(X\bigr)$.
\end{lem}
\begin{proof}
For $X=f(z)g(Y)$ with $f$ and $g$ Laurent polynomials, we have
\[\hbox{Ev}_{\underline{\tilde{t}},q}(\nu_{\underline{t},q}(X))=
\bigl(g(\widetilde{z}^{-1})f(\widetilde{Y}^{-1})(1)\bigr)(k_1^{-1}k_0^{-1})=
f(k_1^{-1}u_1^{-1})g(k_1k_0)=\hbox{Ev}_{\underline{t},q}(X),\]
and for $X=f(z)T_1g(Y)$ we have
\[\hbox{Ev}_{\underline{\tilde{t}},q}(\nu_{\underline{t},q}(X))=
f(k_1^{-1}u_1^{-1})k_1g(k_1k_0)=\hbox{Ev}_{\underline{t},q}(X).
\]
Combined with proposition 8.2 we obtain the desired result.
\end{proof}

\subsection {}

%10.5
We associate with the evaluation mappings $\hbox{Ev}_{\underline{t},q}$ and
$\hbox{Ev}_{\underline{\tilde{t}},q}$ two bilinear forms
\[B: \mathcal{H}(S;\underline{t};q)\times \mathcal{H}(S;
\underline{\tilde{t}};q)\rightarrow {\mathbb{C}},\qquad
\widetilde{B}: \mathcal{H}(S;\underline{\tilde{t}};q)\times
\mathcal{H}(S;\underline{t};q)\rightarrow {\mathbb{C}},
\]
which are defined by $B(X,\widetilde{X})=\hbox{Ev}_{\underline{t},q}
\bigl(\nu_{\underline{\tilde{t}},q}(\widetilde{X})X\bigr)$ and
$\widetilde{B}\bigl(\widetilde{X},X\bigr)=
\hbox{Ev}_{\underline{\tilde{t}},q}\bigl(
\nu_{\underline{t},q}(X)\widetilde{X}\bigr)$
for $X\in \mathcal{H}(S;\underline{t};q)$ and
$\widetilde{X}\in \mathcal{H}(S;\underline{\tilde{t}};q)$.

\begin{lem} Let $X,X_1,X_2\in \mathcal{H}(S;\underline{t};q)$ and
$\widetilde{X},\widetilde{X}_1, \widetilde{X}_2\in
\mathcal{H}(S;\underline{\tilde{t}};q)$. Let $f\in {\mathcal{A}}$.

{\bf (i)} $B(X,\widetilde{X})=\widetilde{B}(\widetilde{X},X)$.

{\bf (ii)} $B(X_1X_2,\widetilde{X})=
B\bigl(X_2, \nu_{\underline{t},q}(X_1)\widetilde{X}\bigr)$,
and $B(X,\widetilde{X}_1\widetilde{X}_2)=B\bigl(
\nu_{\underline{\tilde{t}},q}(\widetilde{X}_1)X,\widetilde{X}_2\bigr)$.

{\bf (iii)} $B\bigl((X(f))(z),\widetilde{X}\bigr)=
B\bigl(X.f(z),\widetilde{X}\bigr)$ and
$B\bigl(X, (\widetilde{X}(f))(\widetilde{z})\bigr)=
B(X,\widetilde{X}.f(\widetilde{z})\bigr)$.

{\bf (iv)} $B\bigl(XT_i, \widetilde{X}\bigr)=k_iB(X,\widetilde{X}\bigr)$
for
$i=0,1$.
\end{lem}
\begin{proof}
{\bf (i)} This follows from lemma 10.4 and the fact that
$\nu_{\underline{t},q}$
is an anti-algebra homomorphism with inverse
$\nu_{\underline{\tilde{t}},q}$, see proposition 8.5.

{\bf (ii)} This is an immediate consequence of proposition 8.5.

{\bf (iii)} The first equality is a direct consequence of the
identity $(X(f))(z)(1)=X(f)=\bigl(X.f(z)\bigr)(1)$ in ${\mathcal{A}}$.
The second identity follows from the first and from {\bf (i)}.

{\bf (iv)} By the explicit expressions 2.19 for
the difference-reflection operators
$T_i$, we have $T_i(1)=k_i1$ for $i=0$ and $i=1$.
The identities are now immediate.
\end{proof}

\subsection {}

%10.6
We write $x_m=k_1^{\epsilon(m)}u_1^{\epsilon(m)}q^m$ ($m\in {\mathbb{Z}}$)
for the eigenvalues of the Cherednik-Dunkl operator
$\widetilde{Y}\in {\mathcal{H}}(S;\underline{\tilde{t}};q)$,
see 3.5. We assume from now on that the parameters $(\underline{t},q)$
are such that
$P_m(x_0^{-1};\underline{t};q)=
\hbox{Ev}_{\underline{t},q}(P_m(\cdot;\underline{t};q))\not=0$
and such that $P_m(\gamma_0^{-1};\underline{\tilde{t}};q)
=\hbox{Ev}_{\underline{\tilde{t}},q}(P_m(\cdot;\underline{\tilde{t}};q))
\not=0$ for all $m\in {\mathbb{Z}}$, and similarly for $P_m^+$
($m\in {\mathbb{Z}}_+$). By proposition 10.2 and 10.3, the
corresponding generic conditions on the parameters
can be specified explicitly.
We write $s(\gamma):=\gamma+\gamma^{-1}$ for all $\gamma\in
{\mathbb{C}}^*$.

\begin{Def}
{\bf (i)}
The renormalized non-symmetric Askey-Wilson polynomials are defined by
\[ E_{\gamma_m}(x;\underline{t};q):=\frac{P_m(x;\underline{t};q)}
{P_m(x_0^{-1};\underline{t};q)},\qquad m\in {\mathbb{Z}}.
\]
In other words, the non-symmetric Askey-Wilson are normalized such that
they take the
value one at $x=x_0^{-1}$.

{\bf (ii)} The renormalized symmetric Askey-Wilson polynomials are
defined by
\[ E_{s(\gamma_m)}^+\bigl(x;\underline{t};q\bigr):=
\frac{P_m^+(x;\underline{t};q)}{P_m^+(x_0;\underline{t};q)},\qquad
m\in {\mathbb{Z}}_+.
\]
In other words, the symmetric Askey-Wilson polynomials are
normalized such that they take the value one at $x=x_0^{\pm 1}$.
\end{Def}
Observe
that $C_+E_{\gamma_m}=E_{s(\gamma_m)}^+$ for $m\in {\mathbb{Z}}$,
where $C_+=(1+k_1^2)^{-1}(1+k_1T_1)$ is the
idempotent defined in 5.1, since $(T_1f)(k_1^{-1}u_1^{-1})=k_1$
for all $f\in {\mathcal{A}}$.

\subsection {}

%10.7
In the following theorem we prove the duality between the
geometric parameter $x=x_n$ and the spectral parameter $\gamma=\gamma_m$
for the renormalized (non-)sym\-me\-tric Askey-Wilson polynomials, see
Sahi \cite[section 7]{S} for the result in the higher rank setting.

\begin{thm}[Duality]
{\bf (i)} For all $m,n\in {\mathbb{Z}}$ and $f\in {\mathcal{A}}$, we have
\[
f(\gamma_m^{-1})=\widetilde{B}\bigl(f(\widetilde{z}),
E_{\gamma_m}(z;\underline{t};q)\bigr),
\qquad
f(x_n^{-1})=B\bigl(f(z),
E_{x_n}(\widetilde{z};\underline{\tilde{t}};q)\bigr).
\]
In particular, $E_{\gamma_m}(x_n^{-1};\underline{t};q)=
E_{x_n}(\gamma_m^{-1};\underline{\tilde{t}};q)$
for all $m,n\in {\mathbb{Z}}$.

{\bf (ii)} For all $m,n\in {\mathbb{Z}}$ and $f\in {\mathcal{A}}^W$,
we have
\[f(\gamma_m)=\widetilde{B}\bigl(f(\widetilde{z}),
E_{s(\gamma_m)}^+\bigl(z;\underline{t};q\bigr)\bigr),
\qquad
f(x_n)=B\bigl(f(z),
E_{s(x_n)}^+\bigl(\widetilde{z};\underline{\tilde{t}};q\bigr)
\bigr).
\]
In particular, $E_{s(\gamma_m)}^+\bigl(x_n;\underline{t};q\bigr)=
E_{s(x_n)}^+\bigl(\gamma_m;\underline{\tilde{t}};q\bigr)$
for all $m,n\in {\mathbb{Z}}$.
\end{thm}
\begin{proof}
{\bf (i)}
The second statement follows from the first by taking
$f=E_{x_n}(\cdot;\underline{\tilde{t}};q)$ in the first equality and
$f=E_{\gamma_m}(\cdot;\underline{t};q)$ in the second equality
and using lemma 10.5{\bf (i)}.

For the first equality, observe that
\begin{equation*}
\begin{split}
\widetilde{B}\bigl(f(\widetilde{z}), E_{\gamma_m}(z)\bigr)&=
\widetilde{B}\bigl(1,f(Y^{-1})E_{\gamma_m}(z)\bigr)=
\widetilde{B}\bigl(1,(f(Y^{-1})E_{\gamma_m})(z)\bigr)\\
&=f(\gamma_m^{-1})\widetilde{B}(1, E_{\gamma_m}(z)\bigr)=
f(\gamma_m^{-1})E_{\gamma_m}(x_0^{-1})=f(\gamma_m^{-1})
\end{split}
\end{equation*}
by application of lemma 10.5{\bf (i)}, {\bf (ii)} and {\bf (iii)}.
The second equality is proved in a similar manner.

{\bf (ii)} The proof is similar to the proof of {\bf (i)}, taking
account of the fact that $f(Y)E_{s(\gamma)}^+=f(\gamma)E_{s(\gamma)}^+$
for all $f(Y)\in {\mathbb{C}}[Y+Y^{-1}]$ by theorem 4.4 and proposition
5.3.
\end{proof}

\subsection {}

%10.8
We write $\sigma=\{\gamma_m\}_{m\in {\mathbb{Z}}}$ for the spectrum
of the Dunkl operator
$Y\in {\mathcal{H}}(S;\underline{t};q)$, see 3.5.
We define an action of the affine Weyl group ${\mathcal{W}}$ on $\sigma$ by
$s_0(\gamma_m)=\gamma_{-m-1}$ and $s_1(\gamma_m)=\gamma_{-m}$
for all $m\in {\mathbb{Z}}$.

The duality between the geometric and spectral parameter of the
renormalized
non-symmetric Askey-Wilson polynomials can be used to explicitly
compute $X(E_{\gamma})$ for $X\in \mathcal{H}(S;\underline{t};q)$
as linear combination of non-symmetric Askey-Wilson polynomials,
cf. 4.1 and 4.2
for $X=T_1$. Observe that $Y(E_{\gamma})=\gamma E_{\gamma}$
and that $T_0^\vee=q^{-1/2}Y^{-1}(T_1^\vee)^{-1}$ by theorem 2.22 and 8.4.
Hence it suffices to expand $T_1(E_{\gamma})$ and
$T_1^\vee(E_{\gamma})$ as linear combination of renormalized
non-symmetric Askey-Wilson polynomials.
The result is as follows.

\begin{prop}
Let $\gamma\in \sigma$, then
\begin{equation*}
\begin{split}
T_1(E_{\gamma})&=k_1E_{\gamma}+k_1^{-1}\frac{(1-k_0k_1\gamma^{-1})
(1+k_0^{-1}k_1\gamma^{-1})}{(1-\gamma^{-2})}\bigl(E_{s_1\gamma}-
E_{\gamma}\bigr),\\
T_1^\vee(E_{\gamma})&=u_1E_{\gamma}+u_1^{-1}\frac{(1-u_0u_1q^{1/2}\gamma)
(1+u_0^{-1}u_1q^{1/2}\gamma)}{(1-q\gamma^2)}
\bigl(E_{s_0\gamma}-E_{\gamma}\bigr).\\
\end{split}
\end{equation*}
\end{prop}
\begin{proof}
The first formula is obviously correct for $\gamma=\gamma_0$.
Let $\gamma_0\not=\gamma\in \sigma$.
By theorem 10.7 and lemma 10.5{\bf (ii)} and {\bf (iii)} we have
\[ (T_1E_{\gamma})(x_m^{-1})=B\bigl(E_{\gamma}(z),
\widetilde{T}_1.E_{x_m}(\widetilde{z};\underline{\tilde{t}};q)\bigr).
\]
By the commutation relation 2.21 between $T_1$ and $f(z)$
($f\in {\mathcal{A}}$) and by
the identity $B(X,\widetilde{X}\widetilde{T}_1)=k_1B(X,\widetilde{X})$ (see
lemma
10.5{\bf (i)} and {\bf (iv)}), we obtain
\begin{equation*}
\begin{split}
(&T_1E_{\gamma})(x_m^{-1})=
k_1B\bigl(E_{\gamma}(z),E_{x_m}(\widetilde{z}^{-1})\bigr)\\
&+\frac{(k_1-k_1^{-1})+(k_0-k_0^{-1})\gamma^{-1}}
{\bigl(1-\gamma^{-2}\bigr)}\bigl(B\bigl(E_{\gamma}(z),
E_{x_m}(\widetilde{z})\bigr)-
B\bigl(E_{\gamma}(z), E_{x_m}(\widetilde{z}^{-1})\bigr)\bigr).
\end{split}
\end{equation*}
Here we have used lemma 10.5{\bf (ii)}, as well as the short-hand notation
$E_{x_m}(\widetilde{z})=E_{x_m}(\widetilde{z};\underline{\tilde{t}};q)$.
By lemma 10.5{\bf (ii)}, {\bf (iii)} and theorem 10.7, we have
\[B\bigl(E_{\gamma}(z), E_{x_m}(\widetilde{z}^{-1})\bigr)=
E_{x_m}(\gamma;\underline{\tilde{t}};q)=E_{\gamma^{-1}}
(x_m^{-1};\underline{t};q)
\]
since $\gamma\not=\gamma_0$. By theorem 10.7 it follows that
the formula for $(T_1E_{\gamma})(x)$ as stated in the proposition is
correct
for $x=x_m^{-1}$ ($m\in {\mathbb{Z}}$), hence it is correct as identity in
${\mathcal{A}}$.

For the second formula we proceed in a similar manner.
First of all, observe that $(T_1^\vee E_{\gamma})(x_m^{-1})=
B\bigl(E_{\gamma}(z), \widetilde{T}_0.E_{x_m}(\widetilde{z};
\underline{\tilde{t}};q)\bigr)$. We use now the commutation
relation 2.21 between $T_0$ and $f(z)$ ($f\in {\mathcal{A}}$) and the
identity
$B\bigl(X,\widetilde{X}\widetilde{T}_0\bigr)=u_1B\bigl(X,\widetilde{X})$,
which follows
from lemma 10.5{\bf (i)} and {\bf (iv)}. Then
\begin{equation*}
\begin{split}
(&T_1^\vee E_{\gamma})(x_m^{-1})=u_1B\bigl(E_{\gamma}(z),
E_{x_m}(q\widetilde{z}^{-1})\bigr)\\
&+\frac{(u_1-u_1^{-1})+(u_0-u_0^{-1})q^{1/2}\gamma}
{(1-q\gamma^2)}\bigl(B\bigl(E_{\gamma}(z), E_{x_m}(\widetilde{z})\bigr)-
B\bigl(E_{\gamma}(z), E_{x_m}(q\widetilde{z}^{-1})\bigr).
\end{split}
\end{equation*}
By lemma 9.2{\bf (ii)}, {\bf (iii)} and theorem 10.7, we have
\[
B\bigl(E_{\gamma}(z), E_{x_m}(q\widetilde{z}^{-1})\bigr)=E_{x_m}(q\gamma;
\underline{\tilde{t}};q)=E_{s_0\gamma}(x_m^{-1};\underline{t};q)
\]
since $q\gamma_n=\gamma_{-n-1}^{-1}$ for all $n\in {\mathbb{Z}}$.
It follows then by direct computations that the formula for
$(T_1^\vee E_{\gamma})(x)$
as stated in the proposition is correct for $x=x_m^{-1}$ ($m\in
{\mathbb{Z}})$,
hence it is correct as identity in ${\mathcal{A}}$.
\end{proof}

\subsection {}

%10.9
The duality for the
(non-)symmetric Askey-Wilson polynomials (see theorem 10.7) can also be
used to
derive recurrence relations for the (non-)symmetric Askey-Wilson
polynomials from the difference(-reflection) equations
$LE_{s(\gamma)}^+=s(\gamma)E_{s(\gamma)}^+$
(respectively  $YE_{\gamma}=\gamma E_{\gamma}$).
The Askey-Wilson $q$-difference equation
$LE_{s(\gamma)}=s(\gamma)E_{\gamma}$
then gives the three term recurrence relation
 \cite[(1.24)--(1.27)]{AW} for the symmetric Askey-Wilson
polynomials (see van Diejen \cite[section 4]{vD}
for the argument in the higher rank setting).

%%%%%%%%%%%%%%%%%%%%%%%%%%%%%%%%%%%%%%%%%%%%%%%%%%%%%%%%%%%%%%%%%%%
%%                                                               %%
%%     The non-symmetric Askey-Wilson transform and its inverse  %%
%%                                                               %%
%%%%%%%%%%%%%%%%%%%%%%%%%%%%%%%%%%%%%%%%%%%%%%%%%%%%%%%%%%%%%%%%%%%

\section{The non-symmetric Askey-Wilson transform and its inverse}

\subsection {}

%11.1
We assume in this section that the parameters
$(\underline{t},q)$ and $(\underline{\tilde{t}},q)$ satisfy the
additional conditions 6.1.

Let $\sigma^\prime$ be the spectrum of $Y^\prime$, so $\sigma^\prime=
\{\gamma_m^\prime\}_{m\in
{\mathbb{Z}}}$ with $\gamma_m^\prime=\gamma_m^{-1}$
for all $m\in {\mathbb{Z}}$, see 3.5.
Let $F=F_{k_0,k_1,q}$ be the functions $g:\sigma^\prime\rightarrow
{\mathbb{C}}$
with finite support. By the non-degeneracy of the bilinear form $\langle
..,. \rangle_{\underline{t},q}$ on ${\mathcal{A}}$ and by the
bi-orthogonality relations 6.7 for the non-symmetric Askey-Wilson
polynomials, we have a bijective linear map
${\mathcal{F}}={\mathcal{F}}_{\underline{t},q}: {\mathcal{A}}\rightarrow F$
defined by
\[ \bigl({\mathcal{F}}_{\underline{t},q}(f)\bigr)(\gamma):=\langle
f,E_{\gamma}^\prime(\,\cdot\,)\rangle_{\underline{t},q},\qquad
f\in {\mathcal{A}},\,\,\, \gamma\in \sigma^\prime,
\]
where
$E_{\gamma}^\prime(\,\cdot\,)=
E_{\gamma}(\,\cdot\,;\underline{t}^{-1};q^{-1})$
($\gamma\in\sigma^\prime$) are the renormalized non-symmetric Askey-Wilson
polynomials with respect to inverse parameters.
\begin{Def}
The bijective map ${\mathcal{F}}: {\mathcal{A}}\rightarrow F$ is
called the non-symmetric Askey-Wilson transform.
\end{Def}

\subsection {}

%11.2
Recall the action of ${\mathcal{W}}$ on $\sigma^\prime$ defined by
$s_0\gamma_m^\prime=\gamma_{-m-1}^\prime$ and
$s_1\gamma_m^\prime=\gamma_{-m}^\prime$
for all $m\in {\mathbb{Z}}$. This induces a left action of ${\mathcal{W}}$
on
$F$ by
$(w\,g)(\gamma)=g(w^{-1}\gamma)$ for $w\in {\mathcal{W}}$, $g\in F$ and
$\gamma\in\sigma^\prime$. Let $\widetilde{T}_i$ ($i=0,1$) and
$\widetilde{z}$
be the linear endomorphisms of $F$ defined by
$(\widetilde{z}g)(\gamma)=\gamma
g(\gamma)$,
\begin{equation*}
\begin{split}
\bigl(\widetilde{T}_0g\bigr)(\gamma)&=u_1g(\gamma)+u_1^{-1}
\frac{(1-u_0u_1q^{1/2}\gamma^{-1})(1+u_0^{-1}u_1q^{1/2}\gamma^{-1})}
{(1-q\gamma^{-2})}\bigl((s_0g)(\gamma)-g(\gamma)\bigr),\\
\bigl(\widetilde{T}_1g\bigr)(\gamma)&=k_1g(\gamma)+k_1^{-1}\frac{(1-k_0k_1
\gamma)
(1+k_0^{-1}k_1\gamma)}{(1-\gamma^2)}\bigl((s_1g)(\gamma)-g(\gamma)\bigr)
\end{split}
\end{equation*}
for all $g\in F$ and all $\gamma\in \sigma^\prime$. Observe that
these formulas can be obtained from the standard action of the
generators $\widetilde{T}_i$ ($i=0,1$) and $\widetilde{z}$ of the double
affine
Hecke algebra
$\mathcal{H}(S;\underline{\tilde{t}};q)$ on ${\mathcal{A}}$ (see
2.19 and 2.22) by formally replacing the
${\mathcal{W}}$-module ${\mathcal{A}}$  by the ${\mathcal{W}}$-module $F$.

\begin{prop}
There is a unique action of $\mathcal{H}(S;\underline{\tilde{t}};q)$
on $F$ such that the generators $\widetilde{z}$ and $\widetilde{T}_i$
\textup{(}$i=0,1$\textup{)} of
$\mathcal{H}(S;\underline{\tilde{t}};q)$
act as the linear endomorphisms defined above.
Furthermore,
\[ {\mathcal{F}}(Xf)=\mu(X){\mathcal{F}}(f),\qquad
X\in \mathcal{H}(S;\underline{t};q),\,\,\, f\in
{\mathcal{A}},
\]
where $\mu$ is the duality isomorphism defined in proposition 8.5.
\end{prop}
\begin{proof}
By proposition 8.7 we have
\[
{\mathcal{F}}(Yf)(\gamma)=\langle f,
(Y^\prime)^{-1}E_{\gamma}^\prime\rangle=
\gamma^{-1}({\mathcal{F}}f)(\gamma)=
\bigl(\widetilde{z}^{-1}{\mathcal{F}}(f)\bigr)(\gamma)=
\bigl(\mu(Y){\mathcal{F}}(f)\bigr)(\gamma)
\]
for all $f\in {\mathcal{A}}$ and all $\gamma\in \sigma^\prime$.

Again by proposition 8.7, we have
${\mathcal{F}}(T_1f)(\gamma)=\langle f,
(T_1^\prime)^{-1}E_{\gamma}^\prime\rangle$. Combined with
proposition 10.8,
we derive that
${\mathcal{F}}(T_1f)=\widetilde{T}_1{\mathcal{F}}(f)=
\mu(T_1){\mathcal{F}}(f)$
for all $f\in {\mathcal{A}}$.

In a similar manner, we derive from proposition 8.7 and proposition
10.8 that
${\mathcal{F}}(T_1^\vee f)=\widetilde{T}_0{\mathcal{F}}(f)=
\mu(T_1^\vee){\mathcal{F}}(f)$ for all $f\in {\mathcal{A}}$.
The proposition now follows since
${\mathcal{F}}$ is bijective and $\mathcal{H}(S;\underline{t};q)$
is generated as an algebra by $Y$, $T_1$ and $T_1^\vee$.
\end{proof}

\subsection {}

%11.3
Observe that the weight function $\Delta(\gamma;
\underline{\tilde{t}};q)$
(see 6.2) has simple poles at $\gamma\in\sigma^\prime$.
We define now a linear map
${\mathcal{G}}={\mathcal{G}}_{\underline{t},q}: F\rightarrow {\mathcal{A}}$
by
\[ {\mathcal{G}}_{\underline{t},q}(g)(x)=\sum_{\gamma\in\sigma^\prime}
g(\gamma)E_{\gamma^{-1}}(x;\underline{t};q)
w(\gamma;\underline{\tilde{t}};q),
\qquad g\in F,
\]
where $w(\gamma)=w(\gamma;\underline{\tilde{t}};q)$ is defined by
\[ w(\gamma;\underline{\tilde{t}};q)=
\underset{y=\gamma}{\hbox{Res}}\left(\frac{\Delta(y;\underline{\tilde{t}};
q)}{y}
\right)\hbox{sgn}(\gamma),\qquad \gamma\in\sigma^\prime
\]
and $\hbox{sgn}(\gamma_m^\prime)=\epsilon(m)$ for $m\in {\mathbb{Z}}$
(see 3.3 for the definition of $\epsilon$). Observe that
$w(\gamma;\underline{\tilde{t}};q)=
\alpha(\gamma;k_1,k_0)w_+(\gamma;\underline{\tilde{t}};q)$
for all $\gamma\in\sigma^\prime$, where
$w_+(\gamma)=w_+(\gamma;\underline{\tilde{t}};q)$ is defined by
\[w_+(\gamma;\underline{\tilde{t}};q)=
\underset{y=\gamma}{\hbox{Res}}
\left(\frac{\Delta_+(y;\underline{\tilde{t}};q)}{y}\right)
\hbox{sgn}(\gamma),\qquad \gamma\in\sigma^\prime,
\]
see 6.2. The weight functions $w(\gamma)$ and $w_+(\gamma)$
can be written out explicitly in terms of
$q$-shifted factorials, see \cite[section 2]{AW} or \cite[section 7.5]{GR}
for
$w_+(\gamma)$.

\subsection {}

%11.4
In the following proposition we determine the intertwining properties
of the $\mathcal{H}(S;\underline{\tilde{t}};q)$-action on $F$
under the linear map ${\mathcal{G}}: F\rightarrow {\mathcal{A}}$.

\begin{prop}
We have
\[ {\mathcal{G}}(Xg)=\mu^{-1}(X)\,{\mathcal{G}}(g),\qquad X\in
\mathcal{H}(S;\underline{\tilde{t}};q),\,\,\, g\in F,
\]
where $\mu$ is the duality isomorphism defined in proposition 8.5.
\end{prop}
\begin{proof}
We write $\widetilde{T}_0=u_1+\widetilde{\phi}_0(\,\cdot\,)(s_0-1)$
and $\widetilde{T}_1=k_1+\widetilde{\phi}_1(\,\cdot\,)(s_1-1)$ with
\begin{equation*}
\begin{split}
\widetilde{\phi}_0(\gamma)&=u_1^{-1}\frac{(1-u_0u_1q^{1/2}\gamma^{-1})
(1+u_0^{-1}u_1q^{1/2}\gamma^{-1})}{(1-q\gamma^{-2})},\\
\widetilde{\phi}_1(\gamma)&=k_1^{-1}\frac{(1-k_0k_1\gamma)(1+k_0^{-1}k_1
\gamma)}{(1-\gamma^2)}.
\end{split}
\end{equation*}
The weight function $w(\gamma;\underline{\tilde{t}};q)$ ($\gamma\in
\sigma^\prime$)
satisfies the fundamental relations
\[
\widetilde{\phi}_i(\gamma)w(\gamma;\underline{\tilde{t}};q)=
\widetilde{\phi}_i(s_i\gamma)w(s_i\gamma;\underline{\tilde{t}};q),\qquad
\gamma\in\sigma',\,\,\, i=0,1.
\]
This follows easily from the explicit
expression for the weight function $\Delta$, see 6.2 (compare
also with the proof of proposition 6.6). It follows that
\begin{equation*}
\begin{split}
{\mathcal{G}}(\widetilde{T}_0g)(x)&=
\sum_{\gamma\in\sigma^\prime}(\widetilde{T}_0g)(\gamma)E_{\gamma^{-1}}(x)\,
w(\gamma)\\
&=\sum_{\gamma\in\sigma}g(\gamma^{-1})\bigl(u_1E_{\gamma}(x)
+\widetilde{\phi}_0(\gamma^{-1})(E_{s_0\gamma}(x)-E_{\gamma}(x))\bigr)\,
w(\gamma^{-1})\\
&=\sum_{\gamma\in\sigma}g(\gamma^{-1})\bigl(T_1^\vee
E_{\gamma}\bigr)(x)\,w(\gamma^{-1})=
\bigl(T_1^\vee {\mathcal{G}}(g)\bigr)(x)
\end{split}
\end{equation*}
for all $g\in F$ by proposition 10.8. Similarly, we obtain
${\mathcal{G}}(\widetilde{T}_1g)(x)=\bigl(T_1{\mathcal{G}}(g)\bigr)(x)$
for all $g\in F$ by proposition 10.8. Furthermore, it is immediate that
${\mathcal{G}}(\widetilde{z}g)(x)= \bigl(Y^{-1}({\mathcal{G}}(g)\bigr)(x)$
for all $g\in F$. We conclude that ${\mathcal{G}}(Xg)=
\mu^{-1}(X){\mathcal{G}}(g)$ for all $g\in F$ if $X=\widetilde{z}$
or $X=\widetilde{T}_i$ for $i=0,1$. The proposition follows, since these
elements generate
$\mathcal{H}(S;\underline{\tilde{t}};q)$ as an algebra.
\end{proof}

\subsection {}

%11.5
Combining proposition 11.2 and proposition 11.4 leads to the
following main result of this section.

\begin{thm}
{\bf (i)} We have
${\mathcal{G}}_{\underline{t},q}\circ
{\mathcal{F}}_{\underline{t},q}=c_{\underline{t},q}\hbox{Id}_{\mathcal{A}}$
and
${\mathcal{F}}_{\underline{t},q}\circ
{\mathcal{G}}_{\underline{t},q}=c_{\underline{t},q}\hbox{Id}_{F}$ with the
constant
$c_{\underline{t},q}$ given by
$c_{\underline{t},q}=w(\gamma_0^{-1};\underline{\tilde{t}};q)\,
\langle 1,1\rangle_{\underline{t},q}$.

{\bf (ii)} For all $\gamma\in\sigma$ we have
\[\frac{\langle E_{\gamma},E_{\gamma^{-1}}^\prime\rangle_{\underline{t},q}}
{\langle
1,1\rangle_{\underline{t},q}}=\frac{w(\gamma_0^{-1};
\underline{\tilde{t}};q)}
{w(\gamma^{-1};\underline{\tilde{t}};q)}.
\]
\end{thm}
\begin{proof}

{\bf (i)} Let $f\in {\mathcal{A}}$, then
\[{\mathcal{G}}\bigl({\mathcal{F}}f\bigr)=
{\mathcal{G}}\bigl({\mathcal{F}}(f(z)1)\bigr)
=f(z)\bigl({\mathcal{G}}({\mathcal{F}}(1))\bigr)
\]
by proposition 11.2 and proposition 11.4,
where $1\in {\mathcal{A}}$ is the function identically equal to one.
By the definitions of ${\mathcal{F}}$ and ${\mathcal{G}}$ and by
the
bi-orthogonality relations for the non-symmetric
Askey-Wilson polynomials (see proposition 6.7) we have
${\mathcal{G}}_{\underline{t},q}({\mathcal{F}}_{\underline{t},q}(1))=
c_{\underline{t},q}1$ with the
constant $c_{\underline{t},q}$ as given in the statement of the
theorem. Hence ${\mathcal{G}}\circ
{\mathcal{F}}=c\,\hbox{Id}_{\mathcal{A}}$.
The identity ${\mathcal{F}}\circ {\mathcal{G}}=c\,\hbox{Id}_F$ follows then
immediately from the fact that ${\mathcal{F}}: {\mathcal{A}}\rightarrow F$
is a bijection.

{\bf (ii)} Let $\gamma\in \sigma$.
By {\bf (i)}, we have ${\mathcal{G}}({\mathcal{F}}(E_{\gamma}))=
c\,E_{\gamma}$.
On the other hand, by the explicit definitions of ${\mathcal{F}}$
and ${\mathcal{G}}$ and by the bi-orthogonality relations for the
non-symmetric Askey-Wilson polynomials (see proposition 6.7),
we have
\[
{\mathcal{G}}_{\underline{t},q}({\mathcal{F}}_{\underline{t},q}(E_{\gamma}))
=
w(\gamma^{-1};\underline{\tilde{t}};q)\langle E_{\gamma},
E_{\gamma^{-1}}^\prime\rangle_{\underline{t},q}\,E_{\gamma}.
\]
Comparing coefficients of $E_{\gamma}$ leads to the desired
result.
\end{proof}

\subsection {}

%11.6
Let $F^W\subset F$ be the $W$-invariant functions in $F$, i.e.
the functions $f\in F$ satisfying $f=s_1f$. Equivalently, $F^W$
consists of the functions $f\in F$ satisfying $\widetilde{C}_+\,f=f$,
where $\widetilde{C}_+=(1+k_1^2)^{-1}(1+k_1\widetilde{T}_1)$.

Let ${\mathcal{F}}_+$ be the restriction of the non-symmetric
Askey-Wilson transform ${\mathcal{F}}$ to ${\mathcal{A}}^W\subset
{\mathcal{A}}$
and let ${\mathcal{G}}_+$ be the restriction of ${\mathcal{G}}$ to
$F^W\subset
F$.

\begin{prop}
${\mathcal{F}}_+$ is a linear bijection from ${\mathcal{A}}^W$ to
$F^W$ with inverse $c^{-1}{\mathcal{G}}_+$, where
$c=c_{\underline{t},q}$ is the constant defined in theorem 11.5.
Furthermore, the constant $c$ can be rewritten as
$c=\frac{1}{2}(1+k_1^2)^2\, w_+(\gamma_0^{-1};\underline{\tilde{t}};q)\,
\bigl(1,1\bigr)_{\underline{t},q}$,
and
\begin{equation*}
\begin{split}
{\mathcal{F}}_+(f)(\gamma)&=
\frac{1}{2}(1+k_1^2)\,
\bigl(f,E_{s(\gamma)}^+(\,\cdot;\underline{t};q)\bigr)_{\underline{t},q},\\
{\mathcal{G}}_+(g)(x)&=(1+k_1^2)\sum_{m\in {\mathbb{Z}}_+}
g(\gamma_m^\prime)E_{s(\gamma_m)}^+(x;\underline{t};q)\,w_+(\gamma_m^\prime;
\underline{\tilde{t}};q)
\end{split}
\end{equation*}
for all $f\in {\mathcal{A}}^W$ and all $g\in F^W$.
\end{prop}
\begin{proof}
By lemma 6.4 we have $\langle 1,1\rangle_{\underline{t},q}=
\frac{1}{2}(1+k_1^2)\bigl(1,1\bigr)_{\underline{t},q}$. Furthermore, we
have
$w(\gamma_0^{-1};\underline{\tilde{t}};q)=
(1+k_1^2)w_+(\gamma_0^{-1};\underline{\tilde{t}};q)$ since
$\alpha(\gamma_0^{-1};k_1,k_0)=1+k_1^2$. This gives
the alternative formula for the constant $c_{\underline{t},q}$.

Observe that $E_{s(\gamma)}^+(x;\underline{t}^{-1};q^{-1})=E_{s(\gamma)}^+(x;
\underline{t};q)$ by
proposition 6.9 and by the $W$-invariance of $E_{s(\gamma)}^+$.
By lemma 6.4, 10.6 and the fact
that $C_+^*=C_+^\prime$ (see the proof of proposition 6.8),
we then derive for $f\in {\mathcal{A}}^W$ and $\gamma\in \sigma^\prime$
that
\begin{equation*}
\begin{split}
{\mathcal{F}}_+(f)(\gamma)=\langle f,
E_{\gamma}^\prime\rangle_{\underline{t},q}&=
\langle C_+f, E_{\gamma}^\prime\rangle_{\underline{t},q}=
\langle f,C_+^\prime E_{\gamma}^\prime\rangle_{\underline{t},q}\\
&=\langle f,E_{s(\gamma)}^+\rangle_{\underline{t},q}=
\frac{1}{2}(1+k_1^2)\bigl(f,E_{s(\gamma)}^+\bigr)_{\underline{t},q}.
\end{split}
\end{equation*}
In particular, we have ${\mathcal{F}}({\mathcal{A}}^W)\subset
F^W$. For $0\not=m\in {\mathbb{Z}}$ we have
\[w(\gamma_m^\prime;\underline{\tilde{t}};q)+w(\gamma_{-m}^\prime;
\underline{\tilde{t}};q)
=(1+k_1^2)w_+(\gamma_m^\prime;\underline{\tilde{t}};q).
\]
Indeed, we use here that $w_+(\gamma;\underline{\tilde{t}};q)=
w_+(\gamma^{-1};\underline{\tilde{t}};q)$ by the $W$-invariance of
the weight function $\Delta_+(\,\cdot\,;\underline{\tilde{t}};q)$,
and that $\alpha(\gamma;k_1,k_0)+
\alpha(\gamma^{-1};k_1,k_0)=1+k_1^2$, see 6.4.
Hence we obtain for $g\in F^W$,
\begin{equation*}
\begin{split}
{\mathcal{G}}_+(g)(x)={\mathcal{G}}(\widetilde{C}_+g)(x)&=
\bigl(C_+{\mathcal{G}}(g)\bigr)(x)
=\sum_{\gamma\in\sigma^\prime}g(\gamma)E_{s(\gamma)}^+(x;\underline{t};q)\,
w(\gamma;\underline{\tilde{t}};q)\\
&=(1+k_1^2)\sum_{m\in {\mathbb{Z}}_+}
g(\gamma_m^\prime)E_{s(\gamma_m)}^+(x;\underline{t};q)\,
w_+(\gamma_m^\prime;\underline{\tilde{t}};q).
\end{split}
\end{equation*}
In particular, ${\mathcal{G}}(F^W)\subset {\mathcal{A}}^W$.
Combined with proposition 11.5, this completes the proof of the
proposition.
\end{proof}

\begin{Def}
The bijection ${\mathcal{F}}_+: {\mathcal{A}}^W\rightarrow F^W$ is
called the symmetric Askey-Wilson transform.
\end{Def}

\subsection {}

%11.7
We can repeat now
the proof of theorem 11.5{\bf (ii)} for the symmetric Askey-Wilson
transform ${\mathcal{F}}_+$, using the alternative descriptions for
${\mathcal{F}}_+$ and ${\mathcal{G}}_+$ as given in proposition 11.6.
This gives the following result on the quadratic norms of the
symmetric Askey-Wilson polynomials.

\begin{cor}
For all $\gamma\in\sigma$ we have
\[ \frac{\bigl(E_{s(\gamma)}^+,
E_{s(\gamma)}^+\bigr)_{\underline{t},q}}{\bigl(1,1\bigr)_{\underline{t},q}}
=\frac{w_+(\gamma_0^{-1};\underline{\tilde{t}};q)}{w_+(\gamma^{-1};
\underline{\tilde{t}};q)}.
\]
\end{cor}

%%%%%%%%%%%%%%%%%%%%%%%%%%%%%%%%%%%%%%%%%%%%%%%%%%%%%%%%%%%%%%%%%%%%%
%%                                                                 %%
%%        The fundamental shift operator and the constant term     %%
%%                                                                 %%
%%%%%%%%%%%%%%%%%%%%%%%%%%%%%%%%%%%%%%%%%%%%%%%%%%%%%%%%%%%%%%%%%%%%%

\section{The fundamental shift operator and the constant term}

\subsection {}

%12.1
In theorem 11.5 and corollary 11.7 we have obtained explicit
expressions of $\langle E_{\gamma},
E_{\gamma^{-1}}^\prime\rangle_{\underline{t};q}$
and of $\bigl(E_{s(\gamma)}^+, E_{s(\gamma)}^+\bigr)_{\underline{t},q}$
in terms of the constant term $\langle 1,1\rangle_{\underline{t},q}=
\frac{1}{2}(1+k_1^2)\bigl(1,1\bigr)_{\underline{t},q}$ for all $\gamma\in
\sigma$. The constant term $\bigl(1,1\bigr)_{\underline{t},q}$
is the well-known Askey-Wilson integral, which has been evaluated
in many different ways, see for instance \cite{AW}, \cite{KM}, \cite{Ko2}
and \cite{R}. We give in this section yet another proof for the
evaluation of $\bigl(1,1\bigr)$ using shift operators.

\subsection {}

%12.2
In the following lemma we define explicit
linear maps from symmetric Laurent polynomials to anti-symmetric
Laurent polynomials and conversely in terms of the Cherednik-Dunkl operator
$Y$.
Recall the definition of the idempotents $C_{\pm}\in H_0\subset
H=H(R;k_0,k_1)$,
see 5.1.

\begin{lem}
Let $h_{\pm}(Y)=h_{\pm}(Y;k_0,k_1)\in H(R;k_0,k_1)$ be defined by
\[ h_{\pm}(Y)=Y^{\mp 1}(Y^{\pm 1}-k_0k_1)(Y^{\pm 1}+k_0^{-1}k_1).
\]

{\bf (i)} We have $h_+(Y)C_+=C_-h_+(Y)C_+$ and $h_-(Y)C_-=C_+h_-(Y)C_-$ in
$H$,
i.e.
$h_{\pm}(Y){\mathcal{A}}_{\pm}\subseteq {\mathcal{A}}_{\mp}$
under the action of the fundamental representation $\pi_{\underline{t},q}$.

{\bf (ii)} We have $C_{\pm}h_{\pm}(Y)C_{\mp}=-h_{\mp}(Y)C_{\mp}$ in $H$.
\end{lem}
\begin{proof}
This follows by a straightforward computation
from Lusztig's formula 2.16, together with the fact that
$(T_1\mp k_1^{\pm 1})C_{\pm}=0$ in $H_0\subset H$.
\end{proof}

\subsection {}

%12.3
By lemma 12.2 and lemma 7.1 we have well-defined linear endomorphisms
$G_{\pm}(\underline{t};q): {\mathcal{A}}^W\rightarrow
{\mathcal{A}}^W$
defined by
\[
\bigl(G_+f\bigr)(x)=\delta(x)^{-1}\bigl(h_+(Y)f\bigr)(x),\qquad
\bigl(G_-f\bigr)=\bigl(h_-(Y)(\delta.f)\bigr)(x),\qquad f\in
{\mathcal{A}}^W,
\]
where $h_{\pm}(Y)$ act under the fundamental representation
$\pi_{\underline{t},q}$. Observe that $G_-$ can be realized as the
element $h_-(Y)\delta(z)$ in
$\mathcal{H}(S;\underline{t};q)$.

\begin{prop}
For $m\in {\mathbb{N}}$, we have
\begin{equation*}
\begin{split}
&G_+(\underline{t};q)P_m^+(\cdot;\underline{t};q)=
h_+(\gamma_m;k_0,k_1)P_{m-1}^+(x;k_0,qk_1,u_0,u_1;q),\\
&G_-(\underline{t};q)P_{m-1}^+(\cdot;k_0,qk_1,u_0,u_1;q)=
h_-(\gamma_m;k_0,k_1)P_m^+(x;\underline{t};q).
\end{split}
\end{equation*}
\end{prop}
\begin{proof}
Let $m\in {\mathbb{N}}$.
Since ${\mathcal{A}}(m)=\hbox{span}\{P_m^+, P_m^-\}$ is an
$H$-module
with ${\mathcal{A}}_-(m)={\mathcal{A}}(m)\cap {\mathcal{A}}_-=
\hbox{span}\{P_m^-\}$,
we have $h_+(Y;k_0,k_1)P_m^+(\cdot;\underline{t};q)=
c_+(m)P_m^-(\cdot;\underline{t};q)$ for some constant $c_+(m)$,
see theorem 4.4, proposition 5.3 and lemma 12.2. Comparing leading terms
using proposition 3.4, we see that $c_+(m)=h_+(\gamma_m)$.
By the generalized Weyl character formula, see proposition 7.3, we obtain
\[ G_+(\underline{t};q)P_m^+(\cdot;\underline{t};q)=
h_+(\gamma_m)P_{m-1}^+(\cdot;k_0,qk_1,u_0,u_1;q).
\]
The shift property of $G_-$ is proved in a similar manner.
\end{proof}

\subsection {}

%12.4
In the remainder of this section we assume that the parameters
satisfy the additional conditions 6.1.
We write
\[
\bigl(G_+(\underline{t}^{-1};q^{-1})f\bigr)(x)=\delta^\prime(x)^{-1}
\bigl(h_+(Y^\prime;k_0^{-1},k_1^{-1})f\bigr)(x),\qquad f\in {\mathcal{A}}
\]
and $G_-(\underline{t}^{-1};q^{-1})=h_-(Y^\prime;k_0^{-1},k_1^{-1})
\delta^\prime(z^\prime)\in {\mathcal{H}}(S;\underline{t}^{-1};q^{-1})$
for the shift-operators with respect to inverse parameters, where
$\delta^\prime(x)=\delta(x;k_0^{-1},k_1^{-1})$ (see 7.2).
The two shift operators $G_+$ and $G_-$ are each-others adjoint in
the following sense.
\begin{prop}
For all $f,g\in {\mathcal{A}}^W$, we have
\begin{equation*}
\begin{split}
\bigl(G_-(\underline{t};q)f,g\bigr)_{\underline{t},q}&=
\bigl(f, G_+(\underline{t}^{-1};
q^{-1})g\bigr)_{k_0,qk_1,u_0,u_1,q},\\
\bigl(G_+(\underline{t};q)f,g\bigr)_{k_0,qk_1,u_0,u_1,q}&=
k_1^4\bigl(f,G_-(\underline{t}^{-1};q^{-1})g\bigr)_{\underline{t},q}.
\end{split}
\end{equation*}
\end{prop}
\begin{proof}
Let $f,g\in {\mathcal{A}}^W$. By lemma 6.4, proposition 6.6 and
lemma 7.1 we have
\begin{equation*}
\begin{split}
\frac{1}{2}(1+k_1^2)\bigl(G_-(\underline{t};q)f,&g\bigr)_{\underline{t},q}=
\langle G_-(\underline{t};q)f,g\rangle_{\underline{t},q}\\
&=\langle \delta(z)f,h_-((Y^\prime)^{-1})g\rangle_{\underline{t},q}=
\langle \delta(z)f, C_-^\prime h_-((Y^\prime)^{-1})C_+^\prime
g\rangle_{\underline{t},q},
\end{split}
\end{equation*}
where $C_{\pm}^\prime\in {\mathcal{H}}(S;\underline{t}^{-1};q^{-1})$
are the images of the primitive idempotents of $H_0(k_1^{-1})$ under
$\pi_{\underline{t}^{-1},q^{-1}}$, see 5.1
(compare with the proof of proposition 6.8).
Now observe that
$h_-((Y^\prime)^{-1};k_0,k_1)=-k_1^2h_-(Y^\prime;k_0^{-1},k_1^{-1})$,
hence
\[
C_-^\prime h_-((Y^\prime)^{-1})C_+^\prime=-k_1^2
C_-^\prime h_-(Y^\prime;k_0^{-1},k_1^{-1})C_+^\prime=
k_1^2h_+(Y^\prime;k_0^{-1},k_1^{-1})C_+^\prime
\]
 by lemma 12.2. Consequently, we obtain
\begin{equation*}
\begin{split}
\frac{1}{2}(1+k_1^2)\bigl(G_-(\underline{t};q)f,g\bigr)_{\underline{t},q}&=
k_1^2\langle \delta(z)f,
h_+(Y^\prime;k_0^{-1},k_1^{-1})g\rangle_{\underline{t},q}\\
&=k_1^2\langle \delta(z)f, \delta^\prime(z^\prime)G_+(\underline{t}^{-1};
q^{-1})g\rangle_{\underline{t},q}\\
&=\frac{1}{2}(1+k_1^2)
\bigl(f,G_+(\underline{t}^{-1};q^{-1})g\bigr)_{k_0,qk_1,u_0,u_1,q}
\end{split}
\end{equation*}
where the last equality follows from lemma 7.2.
The second formula is proved in a similar manner.
\end{proof}

\subsection {}

%12.5
We write $\nu(f)=\nu_{a,b,c,d}(f)=\bigl(f,f\bigr)_{\underline{t},q}$
for the ``quadratic norm''
of $f$ with respect to the bilinear form
$\bigl(.,.\bigr)_{\underline{t},q}$,
where
$(a,b,c,d)$ is the reparametrized multiplicity function, see 5.7.
We use the short-hand notation
\[
\nu_{a,b,c,d}(P_m^+)=\nu_{a,b,c,d}(P_m^+(\cdot;a,b,c,d)),\qquad m\in
{\mathbb{Z}}_+.
\]
\begin{cor}
For $m\in {\mathbb{N}}$, we have
\[
\nu_{a,b,c,d}(P_m^+)=\frac{(1-q^m)(1-q^{m-1}cd)}
{(1-q^mab)(1-q^{m-1}abcd)}
\nu_{qa,qb,c,d}(P_{m-1}^+).
\]
\end{cor}
\begin{proof}
This is an immediate consequence of proposition 12.3, proposition 12.4
and proposition 6.9.
\end{proof}

\subsection {}

%12.6
Observe that the Askey-Wilson second order $q$-difference operator $L$
(see 5.8) and its eigenvalues
$\gamma_m+\gamma_m^{-1}$ ($m\in {\mathbb{Z}}_+$)
are symmetric in $a,b,c,d$. It follows that the symmetric
Askey-Wilson polynomials $P_m^+(x;a,b,c,d)$ ($m\in
{\mathbb{Z}}_+$) are symmetric in the four parameters $a,b,c,d$.
Hence corollary 12.5 can be
reformulated with the special role of $(a,b)$ replaced by an arbitrary
pair of the four parameters $a$, $b$, $c$, $d$.
This leads to the following result.
\begin{cor}
Let $k,l,m,n\in {\mathbb{Z}}_+$ and set $t=k+l+m+n$. Then
\begin{equation*}
\begin{split}
&\frac{\nu_{a,b,c,d}(P_t^+)}
{\nu_{q^{2k}a,q^{2l}b,q^{2m}c,q^{2n}d}(1)}=\\
&\,\,\frac{\bigl(q,q^{2t-1}abcd,q^{2k+2l}ab, q^{2k+2m}ac, q^{2k+2n}ad,
q^{2l+2m}bc, q^{2l+2n}bd, q^{2m+2n}cd;q\bigr)_{\infty}}
{\bigl(q^{t+1}, q^{t-1}abcd, q^tab, q^tac, q^tad, q^tbc, q^tbd, q^tcd;
q\bigr)_{\infty}}.
\end{split}
\end{equation*}
\end{cor}
\begin{proof}
We write
\[\frac{\nu_{a,b,c,d}(P_k^+)}
{\nu_{q^2a,b,c,d}(P_{k-1}^+)}=\frac{\nu_{a,b,c,d,q}(P_k^+)}
{\nu_{qa,qb,c,d}(P_{k-1}^+)}
\frac{\nu_{qa,qb,c,d}(P_{k-1}^+)}
{\nu_{qa,b,q^{-1}c,d}(P_{k}^+)}
\frac{\nu_{qa,b,q^{-1}c,d}(P_{k}^+)}
{\nu_{q^2a,b,c,d}(P_{k-1}^+)}
\]
for $k\in {\mathbb{N}}$ and use the symmetry in the parameters $a,b,c,d$
and
corollary 12.5 to obtain
\[\frac{\nu_{a,b,c,d}(P_k^+)}
{\nu_{q^2a,b,c,d}(P_{k-1}^+)}=\frac{(1-q^k)(1-q^{k-1}bc)(1-q^{k-1}bd)
(1-q^{k-1}cd)}
{(1-q^{k-1}abcd)(1-q^kab)(1-q^{k}ac)(1-q^kad)}
\]
for $k\in {\mathbb{N}}$. Now use again the symmetry in the parameters
$a,b,c,d$
and complete induction with respect to $k,l,m$ and $n$ to obtain the
desired result.
\end{proof}

\subsection {}

%12.7
Corollary 12.6 relates the quadratic norm $\nu(P_m^+)$
to the constant term $\nu(1)$, but it can also be used to
evaluate the Askey-Wilson integral $\nu(1)$ itself.
The Askey-Wilson integral was
evaluated for the first time by Askey and Wilson \cite[theorem 2.1]{AW}
(see also e.g. \cite{KM}, \cite{Ko2} and \cite{R} for alternative proofs).

\begin{thm}[Constant term evaluation]
We have
\[
\nu_{a,b,c,d}(1)=\frac{2\bigl(abcd;q\bigr)_{\infty}}
{\bigl(q,ab,ac,ad,bc,bd,cd;q\bigr)_{\infty}}.
\]
\end{thm}
\begin{proof}
Let
$(a,b,c,d)=(1,-1,q^{\frac{1}{2}}, -q^{\frac{1}{2}})$,
then we have $\Delta_+(x)\equiv 1$ for the corresponding weight
function of the bilinear form $\bigl( .,. \bigr)$. Hence
$P_k^+(x)=x^k+x^{-k}$ ($k\in {\mathbb{N}}$) for the corresponding
symmetric Askey-Wilson polynomial. Furthermore,
\[
\nu_{1,-1,q^{\frac{1}{2}}, -q^{\frac{1}{2}}}(P_k)=2,\quad k\in
{\mathbb{N}}.
\]
Combined with corollary 12.6 this implies that the theorem is
correct for all parameter values
$(a,b,c,d)=\bigl(q^{2k}, -q^{2l}, q^{\frac{1}{2}+2m},
-q^{\frac{1}{2}+2n}\bigr)$
with $k,l,m,n\in {\mathbb{Z}}_+$ and $k+l+m+n\in {\mathbb{N}}$
(we use here that the formula in corollary 12.6 extends to this particular
choice of
parameter values by continuity).
The proof is now completed by analytic continuation.
\end{proof}

\subsection {}

%12.8
The constant term evaluation (theorem 12.7) and corollary 12.5
yield an explicit evaluation of $\nu_{a,b,c,d}(P_m^+)$ for $m\in
{\mathbb{Z}}_+$ which is in accordance with Askey and Wilson's
result \cite[theorem 2.2]{AW}. As remarked in 7.5, this then yields
explicit evaluations for all the diagonal terms in the
bi-orthogonality relations of proposition 6.7 and proposition 6.8.

Another way to obtain the diagonal terms explicitly is by using
corollary 11.7 and 10.3 (respectively theorem 11.5 and proposition
10.2) to reduce the diagonal terms for the (non-)symmetric
Askey-Wilson polynomials to the constant term evaluation (theorem 12.7).
Explicitly, we obtain the following formulas for the diagonal
terms:
\[\bigl(P_m^+, P_m^+\bigr)=\frac{2\bigl(q^{2m-1}abcd,
q^{2m}abcd;q\bigr)_{\infty}}
{\bigl(q^{m+1}, q^mab, q^mac, q^mad, q^mbc, q^mbd, q^mcd,
q^{m-1}abcd;q\bigr)_{\infty}}
\]
for $m\in {\mathbb{Z}}_+$,
\[\langle P_m, P_m^\prime\rangle=
\frac{\bigl(q^{2m}abcd, q^{2m}abcd;q\bigr)_{\infty}}
{\bigl(q^{m+1}, q^{m+1}ab, q^mac, q^mad, q^mbc, q^mbd, q^mcd,
q^mabcd;q\bigr)_{\infty}}
\]
for $m\in {\mathbb{Z}}_+$,
\[\langle P_{-m}, P_{-m}^\prime\rangle=
\frac{\bigl(q^{2m-1}abcd, q^{2m-1}abcd;q\bigr)_{\infty}}
{\bigl(q^m,q^mab,q^mac,q^mad,q^mbc,q^mbd,q^{m-1}cd,q^{m-1}abcd;
q\bigr)_{\infty}}
\]
for $m\in {\mathbb{N}}$ and finally
\[
\langle P_m^-, P_m^{- \prime}\rangle=
\frac{ab-1}{ab}\frac{\bigl(q^{2m-1}abcd, q^{2m}abcd;q\bigr)_{\infty}}
{\bigl(q^m, q^{m+1}ab, q^mac, q^mad, q^mbc, q^mbd, q^{m-1}cd,
q^mabcd;q\bigr)_{\infty}}
\]
for $m\in {\mathbb{N}}$.

\subsection {}

%12.9
There is yet another way to relate the
diagonal terms of the non-symmetric Askey-Wilson polynomials
to the constant term $\langle 1,1\rangle$. This method is based on
the Rodrigues type formula for the non-symmetric Askey-Wilson
polynomial (see proposition 9.3), which allows us to compute
the diagonal terms by induction with respect to the
degree of the non-symmetric Askey-Wilson
polynomial. For the induction step, one
needs the following two additional properties of the
intertwiners. The first property is
that $S_0^*=q^{-1}S_0^\prime$ and $S_1^*=S_1^\prime$, where
$S_0^\prime, S_1^\prime\in {\mathcal{H}}(S;\underline{t}^{-1};q^{-1})$ are
the intertwiners with respect to inverse parameters, cf.
proposition 8.7. The second property is
\begin{equation*}
\begin{split}
S_0^2&=q^{-1}u_1^2\prod_{\xi=\pm 1}
(1-u_0^{-1}u_1^{-1}q^{\xi/2}Y^{\xi})(1+u_0u_1^{-1}q^{\xi/2}Y^{\xi}),\\
S_1^2&=k_1^2\prod_{\xi=\pm
1}(1-k_0^{-1}k_1^{-1}Y^{\xi})(1+k_0k_1^{-1}Y^{\xi})
\end{split}
\end{equation*}
which are most easily proved in the image of the duality
isomorphism $\mu$, see Sahi \cite[corollary 5.2]{S}
in the higher rank setting. We leave the details to the reader.

%%%%%%%%%%%%%%%%%%%%%%%%%%%%%%%%%%%%%%%%%%%%%%%%%%%%%%%%%%%%%%%%
%%                                                            %%
%%                       Bibliography                         %%
%%                                                            %%
%%%%%%%%%%%%%%%%%%%%%%%%%%%%%%%%%%%%%%%%%%%%%%%%%%%%%%%%%%%%%%%%

\bibliographystyle{amsplain}

\end{document}